\def\VersionDateTime{1/November/2016, 18:36 GMT+9:00. Version $1.5$}
\documentclass[12pt, leqno]{amsart}
\usepackage{amssymb,amsthm,amsmath,amscd,rsfs}
\usepackage{latexsym}

\newcommand{\NN}{{\mathbb{N}}}
\newcommand{\ZZ}{{\mathbb{Z}}}
\newcommand{\QQ}{{\mathbb{Q}}}
\newcommand{\RR}{{\mathbb{R}}}

\newcommand{\PP}{{\mathbb{P}}}

\newcommand{\OO}{{\mathcal{O}}}

\newcommand{\codim}{\operatorname{codim}}
\newcommand{\ch}{\operatorname{char}}

\newcommand{\Ker}{\operatorname{Ker}}

\newcommand{\Supp}{\operatorname{Supp}}

\newcommand{\Spec}{\operatorname{Spec}}

\newcommand{\cherncl}{{c}}

\newcommand{\zero}{\operatorname{div}}

\newcommand{\Proof}{{\sl Proof.}\quad}
\newcommand{\QED}{{\unskip\nobreak\hfil\penalty50\quad\null\nobreak\hfil
{$\Box$}\parfillskip0pt\finalhyphendemerits0\par\medskip}}
\newcommand{\rest}[2]{\left.{#1}\right\vert_{{#2}}}

\newcommand{\ord}{\operatorname{ord}}

\newcommand{\Tr}{\operatorname{Tr}}

\newcommand{\pr}{\operatorname{pr}}

\title
[Non-density of small points on divisors
and the Bogomolov conjecture
]{Non-density of small points on divisors on abelian varieties
and the Bogomolov conjecture}
\author{Kazuhiko Yamaki}
\date{\VersionDateTime}
\address
{Institute for Liberal Arts and Sciences,
Kyoto University, Kyoto, 606-8501, Japan}
\email{yamaki.kazuhiko.6r@kyoto-u.ac.jp}

\begin{document}

\theoremstyle{plain}
\newtheorem{Theorem}{Theorem}[section]
\newtheorem{Lemma}[Theorem]{Lemma}
\newtheorem{Proposition}[Theorem]{Proposition}
\newtheorem{Corollary}[Theorem]{Corollary}
\newtheorem{Main-Theorem}[Theorem]{Main Theorem}
\newtheorem{Theorem-Definition}[Theorem]{Theorem-Definition}
\theoremstyle{definition}
\newtheorem{Definition}[Theorem]{Definition}
\newtheorem{Remark}[Theorem]{Remark}
\newtheorem{Conjecture}[Theorem]{Conjecture}
\newtheorem{Claim}{Claim}
\newtheorem{Example}[Theorem]{Example}
\newtheorem{Key Fact}[Theorem]{Key Fact}
\newtheorem{ack}{Acknowledgments}       \renewcommand{\theack}{}
\newtheorem*{n-c}{Notation and convention}      
\newtheorem{citeTheorem}[Theorem]{Theorem}
\newtheorem{citeProposition}[Theorem]{Proposition}

\newtheorem{Step}{Step}

\renewcommand{\theTheorem}{\arabic{section}.\arabic{Theorem}}
\renewcommand{\theClaim}{\arabic{section}.\arabic{Theorem}.\arabic{Claim}}
\renewcommand{\theequation}{\arabic{section}.\arabic{Theorem}.\arabic{Claim}}

\def\Pf{\trivlist\item[\hskip\labelsep\textit{Proof.}]}
\def\endPf{\strut\hfill\framebox(6,6){}\endtrivlist}

\def\Pfo{\trivlist\item[\hskip\labelsep\textit{Proof of Proposition~\ref{ch-of-hyp}.}]}
\def\endPfo{\strut\hfill\framebox(6,6){}\endtrivlist}

\maketitle


\begin{abstract}
The Bogomolov conjecture for a curve 
claims finiteness of algebraic points
on the curve
which are small with respect to the canonical height.
Ullmo has established this conjecture over number fields,
and Moriwaki generalized 
it
to the assertion
over 
finitely generated fields over $\QQ$
with respect to arithmetic heights.
As for the case of function fields
with respect to the geometric heights,
Cinkir has proved the conjecture
over function fields of characteristic $0$ 
and of
transcendence degree $1$.
However, the conjecture has been open over
other function fields.

In this paper,
we prove
that
the Bogomolov conjecture
for curves
holds
over any function field.
In fact, we show that 
any non-special closed subvariety of \emph{dimension} $1$ 
in an abelian variety over function fields
has only a finite number of small points.
This result is
a consequence of the investigation of non-density of small points
of closed subvarieties of abelian varieties of \emph{codimension} $1$.
As a by-product,
we remark that the geometric Bogomolov conjecture,
which is a generalization of the Bogomolov conjecture for curves
over function fields,
holds for any abelian variety
of dimension at most $3$.
Combining this result with our previous works,
we see that the geometric Bogomolov conjecture
holds for all abelian varieties
for which the difference between
its nowhere degeneracy rank
and the dimension of its trace is not greater than $3$.
\end{abstract}

\section{Introduction}

\subsection{Bogomolov conjecture for curves} \label{subsect;mainresults}

Let $K$ be a function field or a number field.
Here, a function field means the function field of a normal projective
variety of positive dimension
over an algebraically closed field $k$.
We fix an algebraic closure $\overline{K}$.
We consider the heights on projective varieties over $\overline{K}$
(cf. \S~\ref{sect:NC}).

Let $C$ be a smooth projective curve of genus $g \geq 2$
over $\overline{K}$
and
let $J_C$ be the Jacobian variety of $C$.
Fix a divisor $D$ on $C$ of degree $1$
and let
$
\jmath_{D} : C \hookrightarrow 
J_C
$
be the embedding
defined
by $\jmath_D (x) := x - D$.
For any $P \in J_{C} \left( \overline{K} \right)$
and
any $\epsilon \geq 0$,
set
\[
B_{C} ( P , \epsilon )
:=
\left\{ 
\left. x \in C \left(
\overline{K}
\right) 
\ 
\right|
|| \jmath_D(x) - P ||_{NT} \leq \epsilon
\right\}
,
\]
where $|| \cdot ||_{NT}$ is the semi-norm arising from the N\'eron--Tate
height
 on $J_{C}$
associated to a symmetric theta divisor.
In \cite{bogomolov},
Bogomolov conjectured
the following.
Here, 
in the case where $K$ is a function field,
a curve $C$ over $\overline{K}$
is said to be \emph{isotrivial} if it can be defined over the constant
field $k$.

\begin{Conjecture} \label{conj:BCFC}
Assume that $C$ is non-isotrivial when $K$ is a function field.
Then, for any $P \in J_{C} \left(
\overline{K} \right)$, there exists $r > 0$ such that
$B_{C} ( P , r )$ is
a finite set.
\end{Conjecture}

When $K$ is a number field,
Ullmo proved 
in 1998
that the conjecture holds 
(cf. \cite[Th\'eor\`eme~1.1]{ullmo}).
In the proof, he used the Szpiro--Ullmo--Zhang equidistribution theorem
of small points over an archimedean place.

Suppose that
$K$ is the function field
of a normal projective variety $\mathfrak{B}$
over $k$.
Then
the conjecture has been proved for some $K$ but not for all $K$ yet.
If $\dim (\mathfrak{B}) = 1$ and $\ch (k) = 0$,
after partial answers such as
\cite{moriwakiRBI, faber},
Cinkir
established the conjecture in \cite{cinkir}
by using
\cite{zhang3}. 
For more general $K$,
there are some partial answers such as
\cite{zhang1, moriwaki0, yamaki1, yamaki4}
(under the assumption of $\dim (\mathfrak{B}) = 1$)
and
\cite{gubler2, yamaki6, yamaki7}
(special case of 
the geometric Bogomolov conjecture for abelian varieties),
but 
we have not had a complete answer to the conjecture yet
except for $\dim ( \mathfrak{B} ) = 1$
and $\ch (k) = 0$.

In this paper,
we show that the Bogomolov conjecture holds
over any function field $K$:

\begin{Theorem} [Theorem~\ref{thm:BCFC}] \label{thm:main1intro}
Let $K$ be a function field over any algebraically closed constant field
$k$.
Then Conjecture~\ref{conj:BCFC} holds.
\end{Theorem}


As an application,
we will give an alternative proof of 
the Manin--Mumford
conjecture for curves in positive characteristic
(cf. \S~\ref{subsect:scanlontheorem}).

Note however that Theorem~\ref{thm:main1intro} is not effective 
in contrast to Cinkir's theorem
(when
$K$ is a function field of transcendence degree $1$
over $k$ of characteristic $0$);
Cinkir gave explicitly a positive number $r$ such that
$B_{C} ( P , r )$ is
finite for any $P \in J_C \left( \overline{K} \right)$
when $C$ has semistable reduction over $K$.

\subsection{Non-density of small points of subvarieties of 
abelian varieties}

We should recall the reason why
Cinkir's proof in characteristic $0$
does not work in positive characteristic.
His proof is based on Zhang's work \cite{zhang3}.
Zhang proved there that
if some constants,
called the $\varphi$-invariants,
arising from the reduction graphs
of the curve $C$
are positive,
and if the height of some cycle of 
the triple product $C^{3}$,
called the
Gross--Schoen cycle, are non-negative,
then
the Bogomolov conjecture holds for $C$.
Then
Cinkir 
proved in \cite{cinkir}
the positivity of
$\varphi$-invariants.
This suffices for the Bogomolov conjecture 
in characteristic $0$
because in this case
it is deduced from the Hodge index theorem that
the height of the Gross--Schoen cycle is non-negative.
However in positive characteristic,
the Hodge index theorem,
which is a part of the standard conjecture,
is not known.
Thus
the Bogomolov conjecture in positive characteristic
cannot be deduced in the same way.

The proof of Theorem~\ref{thm:main1intro}
is quite different from 
that of Cinkir's theorem;
we do not use reduction graphs of curves
or the Gross--Schoen cycles.
In our argument,
Theorem~\ref{thm:main1intro}
is a direct consequence of
a more general result 
concerning the non-density of small points on curves in
abelian varieties,
and
this result is deduced from the non-density
of small points on divisors on abelian varieties,
as we will now explain.

Let $A$ be an abelian variety over $\overline{K}$
and let
$L$ be an even ample line bundle on $A$.
Here, ``even'' means the pull-back $[-1]^{\ast} (L)$
by the $(-1)$-times automorphism coincides with $L$.
Let $X$ be a closed subvariety of $A$
and put
\[
X ( \epsilon ; L)
:=
\left\{
x \in X \left(
\overline{K} \right)
\left|
\ 
\widehat{h}_{L}
(x) \leq \epsilon
\right.
\right\}
.
\]
Then
we say that \emph{$X$ has dense small points}
if $X ( \epsilon ; L)$
is dense in $X$ for any $\epsilon > 0$.
The notion of density of small points
is known to be independent of
the choice of such $L$'s
(cf. \cite[Lemma~2.1 and Definition~2.2]{yamaki5}).

Let $K$ be a function field
and
let $A$ be an abelian variety over $\overline{K}$.
We recall the notions of $\overline{K}/k$-trace
and special subvarieties of $A$.
A \emph{$\overline{K}/k$-trace} of $A$ is a pair
$\left( \widetilde{A}^{\overline{K}/k},
\Tr_A \right)$ of an abelian variety over $k$
and a homomorphism $\Tr_A : \widetilde{A}^{\overline{K}/k} \otimes_{k}
\overline{K}
\to A$ such that $\Tr_A$ is the universal one
among the homomorphisms to $A$ from abelian varieties
which can be defined over the constant field $k$
(See \S~\ref{sect:NC} for more details).
It is known that there exists a unique $\overline{K}/k$-trace
of $A$.
Using the $\overline{K}/k$-trace, we defined
the notion of special subvarieties in \cite[\S~2]{yamaki5};
a closed subvariety $X$ of $A$ is said to be \emph{special}
if there exist a closed subvariety $\widetilde{Y}$ of 
$\widetilde{A}^{\overline{K}/k}$,
an abelian subvariety $G$ of $A$, 
and a torsion point $\tau \in A \left(
\overline{K}
\right)$ such that
$X = \Tr_A \left( \widetilde{Y} \otimes_k \overline{K} \right)
+ G + \tau$.

We will prove the following theorem.

\begin{Theorem} [Theorem~\ref{thm:curvecase}]
\label{thm:curvecaseint}
Let $A$ be an abelian variety over $\overline{K}$.
Let $X$ be a closed subvariety of $A$ with $\dim (X) = 1$.
Then if $X$ has dense small points, then it is special.
\end{Theorem}

As is remarked in \cite[\S~8]{yamaki6},
Theorem~\ref{thm:main1intro}
is an immediate consequence of
Theorem~\ref{thm:curvecaseint}.

Theorem~\ref{thm:curvecaseint}
will be deduced
from the following theorem,
which is the crucial result in this paper.

\begin{Theorem} [Theorem~\ref{thm:divisorialcase}]
\label{thm:divisorialcaseint}
Let $A$ be an abelian variety over $\overline{K}$.
Let $X$ be a closed subvariety of $A$
of codimension $1$.
Then if $X$ has dense small points, then it is special.
\end{Theorem}

The basic idea to connect
Theorem~\ref{thm:curvecaseint}
to Theorem~\ref{thm:divisorialcaseint}
is as follows.
Let $X$ be a closed subvariety of $A$ of dimension $1$.
We set $Y_0 := \{ 0 \}$,
and for each positive integer $m$,
let
$Y_m$ denote the sum of $m$ copies of $X - X$.
Then $Y_m$ is an abelian subvariety of $A$ for large $m$,
and let $N$ be the smallest positive integer among such $m$.
Suppose that $X$ has dense small points.
Then we can 
show that there exists 
a torsion
point $\tau \in A \left(
\overline{K}
\right)$
such that
$Y_{N-1}$ or $Y_{N-1} +  ( X - \tau )$
is an effective divisor on $Y_N$,
which we denote by $D$.
Note that 
$D$ 
also has dense small points.
By considering these $D$
and $Y_N$, we can deduce Theorem~\ref{thm:curvecaseint}
from Theorem~\ref{thm:divisorialcaseint}.

These two theorems concern the following conjecture,
called the \emph{geometric Bogomolov conjecture}.
\begin{Conjecture} [Conjecture~2.9 of \cite{yamaki5}] \label{conj:GBCforAVint}
Let $A$ be an abelian variety over $\overline{K}$.
Let $X$ be a closed subvariety of $A$.
Then $X$ has dense small points if and only if 
$X$ is a special subvariety.
\end{Conjecture}

Since 
a special subvariety has dense small points,
the problem is the ``only if'' part.
The geometric Bogomolov conjecture
is the geometric version of the Ullmo--Zhang theorem
\cite{zhang2},
which is also called the Bogomolov conjecture for abelian varieties.
We refer to \S~\ref{subsect:GBC} for more details including some background of
this conjecture.

As a consequence of Theorems~\ref{thm:curvecaseint} and 
\ref{thm:divisorialcaseint},
we see that
the geometric Bogomolov conjecture holds for
$A$ with $\dim (A) \leq 3$
(Corollary~\ref{cor:Aleq3}),
because
any non-trivial closed subvariety
has dimension $1$ or codimension $1$
in this case.
Furthermore,
combining this result
with \cite[Theorem~1.5]{yamaki7}
(cf. Theorem~\ref{thm:yamaki72int}),
we 
see that
the conjecture holds for $A$
with 
$\dim (\mathfrak{m}) \leq  
\dim \left( \widetilde{A}^{\overline{K}/k} \right) + 3$
(Corollary~\ref{cor:mleqt+3} and 
Remark~\ref{rem:trace-finite}),
where $\mathfrak{m}$ is the maximal nowhere degenerate abelian subvariety
of $A$.
See \S~\ref{subsect:GBC} for more details.

\subsection{Idea}

We describe the idea of the proof of 
Theorem~\ref{thm:divisorialcaseint}.
To avoid technical difficulties,
we assume $\dim (\mathfrak{B}) =1$ in this subsection.
Let $A$ be an abelian variety 
of dimension $n$ over $\overline{K}$
and let $X$ be a closed subvariety of $A$ of codimension $1$.
Remark that $X$ has dense small points if and only if $X$ has
canonical height $0$ with respect to an even ample line bundle
(cf. Proposition~\ref{prop:dense-height0}).
Then
what we should show is
that if $X$ is not a special subvariety,
then $X$ has 
positive canonical height.
From Theorem~\ref{thm:yamaki72int},
which is
proved in \cite{yamaki6} and \cite{yamaki7},
together with 
a few arguments,
one sees that
it suffices to show the theorem under the assumption that
$A$ is nowhere degenerate 
and has trivial $\overline{K}/k$-trace,
where ``$A$ is nowhere degenerate'' essentially means
``$A$ has potentially good reduction everywhere''; see
\S~\ref{sect:NC} for the precise terminology.
Further,
noting that any effective divisor is the pullback of an
ample divisor by some homomorphism
(cf. \cite[p.88, Remarks on effective divisors by Nori]{mumford}),
one finds that
it suffices to show
the following assertion
(cf. Proposition~\ref{prop:main1}):
Assume that $A$ is 
a 
nowhere degenerate abelian variety
with trivial $\overline{K}/k$-trace;
let $X$ be a subvariety of $A$ of codimension $1$
and suppose that $X$ is ample;
then $X$ has positive canonical height.

In this description,
we furthermore make the following assumption:
\begin{enumerate}
\renewcommand{\labelenumi}{(\alph{enumi})}
\item
there exists an abelian scheme $f : \mathscr{A} \to \mathfrak{B}$
with zero-section $0_f$
and with
geometric generic fiber $A$,
and
$X$ is defined over $K$;
\item
$0 \notin X$;
\item
$\# k > \aleph_0$,
i.e.,
$k$ has uncountably infinite cardinality.
\end{enumerate}
As we can see later,
those assumptions do not give essential restrictions.

Let us give a sketch of the proof.
Put $D := X + [-1]^{\ast} (X)$, 
where ``$+$'' is the addition of divisors,
and put $L := \OO_A (D)$.
Note that 
$L$ is even and ample.
One sees that
$\widehat{h}_{L} ( D) = 2 \widehat{h}_L (X)$.
Let $\mathscr{D}$ be the closure of $D$ in $\mathscr{A}$.
Since $0 \notin D$,
$0_f^{\ast} ( \mathscr{D} )$ is a well-defined 
effective divisor on $\mathfrak{B}$.
Set $\mathfrak{L} := \OO_{\mathfrak{B}}
\left( 0_f^{\ast} ( \mathscr{D})
\right)$
and
set 
$\mathscr{L} := \OO_{\mathscr{A}} \left( \mathscr{D}
\right) 
\otimes
f^{\ast} \mathfrak{L}^{\otimes (-1)}$.
Then $0_{f}^{\ast} ( \mathscr{L} ) = \OO_{\mathfrak{B}}$,
and
one sees that
\[
\widehat{h}_{L} ( D ) = \deg
\left( 
\cherncl_1 ( \mathscr{L} )^{\cdot n}
\cdot  \mathscr{D} \right)
=
\deg
\left( 
\cherncl_1 ( \mathscr{L} )^{\cdot n}
\cdot  \cherncl_1 \left(
\OO_{\mathscr{A}}
\left(
\mathscr{D}
\right)
\right)
\cdot
\mathscr{A}
\right)
.
\]
Since $A$ has canonical height $0$,
we have $\deg \left( \cherncl_1 ( \mathscr{L} )^{\cdot (n + 1)}
\cdot  \mathscr{A} \right) = 0$
(cf. Remark~\ref{rem:heightabelianvariety}).
It follows that 
\begin{align*}
\deg
\left( 
\cherncl_1 ( \mathscr{L} )^{\cdot n}
\cdot  \cherncl_1 \left(
\OO_{\mathscr{A}}
\left(
\mathscr{D}
\right)
\right)
\cdot
\mathscr{A}
\right)
&=
\deg
\left( 
\cherncl_1 ( \mathscr{L} )^{\cdot (n+1)}
\cdot 
\mathscr{A}
\right)
+
\deg
\left( 
\cherncl_1 ( \mathscr{L} )^{\cdot n}
\cdot 
f^{\ast} \cherncl_1 ( \mathfrak{L})
\cdot
\mathscr{A}
\right)
\\
&=
\deg_L (A) \cdot \deg (\mathfrak{L})
.
\end{align*}
Since $\deg_L (A) > 0$ by the ampleness of $L$,
it remains to show $\deg (\mathfrak{L}) > 0$.
In fact, we see in Lemma~\ref{lem:main}, which is a key lemma,
that $\mathfrak{L}$ is non-trivial.
Since $\mathfrak{L}$ is effective,
we conclude
$\deg (\mathfrak{L}) > 0$.

The outline of the proof of the non-triviality
of $\mathfrak{L}$
is as follows.
We prove the non-triviality by contradiction.
Suppose that it is trivial.
Then one can show that 
there exists a finite covering $\mathfrak{B}' \to \mathfrak{B}$
such that
the complete linear system $| 2 \mathscr{D}' |$
on $\mathscr{A}'$
is base-point free
(cf. Proposition~\ref{prop:basepointfree}),
where $f' : \mathscr{A}' \to \mathfrak{B}'$ and $\mathscr{D}'$
are the base-change of $f$ and $\mathscr{D}$ by this
$\mathfrak{B}' \to \mathfrak{B}$ respectively.
Let
$\varphi : \mathscr{A}' \to Z$
be the surjective morphism
associated to this complete linear system,
where $Z$ is a closed subvariety of the dual space of $|2 \mathscr{D}' |$.
Remark that for any curve $\gamma \subset \mathscr{A}'$,
$\deg \left(
\cherncl_1 ( \mathscr{L}' ) \cdot \gamma
\right) = 0$
if and only if $\varphi ( \gamma )$ is a point,
where $\mathscr{L}'$ is 
the pull-back of $\mathscr{L}$ to $\mathscr{A}'$.
Further, remark that for any $a \in A \left( 
\overline{K} \right)$,
$\widehat{h}_{L} ( a ) = 0$
if and only  if $\deg \left(
\cherncl_1 ( \mathscr{L}' ) \cdot \Delta_{a}
\right) = 0$,
where $\Delta_{a}$ is the closure of $a$
in $\mathscr{A}'$.
Since $\left( \mathscr{L}'\right)^{\otimes 2}$ 
is relatively ample with respect to $f'$,
it follows that $\varphi$ is finite on any fiber,
and 
since the closure of the set of the 
height $0$ points of
$A \left( \overline{K} \right)$
is dense in $\mathscr{A}'$,
we see that $\varphi$ is not generically finite 
on $\mathscr{A}'$.
Since $\# k > \aleph_0$,
it follows that 
$\mathscr{A}'$
contains
uncountably many
irreducible curves which are flat over $\mathfrak{B}'$
and are contracted to a point by $\varphi$.
This means that
$A \left(
\overline{K}
\right)$ has uncountably many
points of height $0$.
On the other hand,
since $A$ has trivial $\overline{K}/k$-trace,
a point of $A \left( 
\overline{K}
\right)$ has height $0$ 
if and only if it is a torsion,
and there are only countably many such points.
This is a contradiction.

\subsection{Organization}
This article consists of six sections including this introduction.
In \S~\ref{sect:preliminary},
we recall canonical heights
and 
their description
in terms of models
when the abelian variety is nowhere degenerate.
In \S~\ref{sect:Bertini},
we prove a version of Bertini's theorem
on normal varieties, which will be used in the
case where $K$ has transcendence degree more than $1$.
In \S~\ref{sect:positivity},
we prove that an effective ample divisor 
on a nowhere degenerate abelian variety with trivial $\overline{K}/k$-trace
has positive canonical height (cf. Proposition~\ref{prop:main1}),
along the idea explained above.
Then, we prove the main results in \S~\ref{sect:GBC}.
In \S~\ref{sect:RGBC}, we apply the results of
\S~\ref{sect:GBC} to obtain partial answers to
the geometric Bogomolov conjecture.
We also remark the relationship between the Manin--Mumford conjecture
over fields of positive characteristic
and the geometric Bogomolov conjecture.

\subsection*{Acknowledgments}
This paper was written in part during my visit
to Regensburg university in March--April 2015,
which was supported by the SFB Higher Invariants.
I thank Professor Walter Gubler for inviting me
and for his hospitality.
I thank him also for a lot of valuable comments for the preliminary draft of this
article.
I thank the referees for valuable comments and suggestions.
This work was partly supported by KAKENHI 26800012.

\section{Preliminary} \label{sect:preliminary}

\subsection{Notation and convention} \label{sect:NC}
In this paper, a natural number means a 
strictly positive integer.
Let $\NN$ denote the set of natural numbers.

Let $F$ be a field 
and let
$X$ be a scheme over $F$.
For a field extension $F'/F$,
we write 
${X} \otimes_F F' := 
{X} \times_{\Spec (F)} \Spec (F')$.
For a morphism $\phi : {X} \to {Y}$ of schemes over $F$,
we write $\phi \otimes_F F' : 
{X} \otimes_F F' \to 
{Y} \otimes_F F'$ for the base-extension
to $F'$.
We call $X$ a
\emph{variety} over $F$ if $X$
is a
geometrically integral scheme 
separated and of finite type over $F$.

Let $f : \mathscr{A} \to S$ be an abelian scheme with zero-section $0_f$.
For any $n \in \ZZ$, let $[n] : \mathscr{A} \to \mathscr{A}$
denote the 
$n$-times 
endomorphism.
Suppose that $S$ is the spectrum of a field,
that is, $f$ is an abelian variety.
Let $L$ be a line bundle on this abelian variety.
We say that $L$ is \emph{even}
if $[-1]^{\ast} (L) \cong L$.
Remark that
in this case,
$[n]^{\ast} (L) \cong L^{\otimes n^{2}}$ holds
by the theorem of the cube
(cf. \cite[\S~6, Corollary~3]{mumford}).

Let $k$ be an algebraically closed field,
$\mathfrak{B}$ a normal projective variety 
of dimension $b \geq 1$ over $k$,
and let $\mathcal{H}$ be an ample line bundle on $\mathfrak{B}$.
Let
$K$ be the function field 
of 
$\mathfrak{B}$
and let $\overline{K}$ be an algebraic closure of $K$.
All of them are fixed throughout this article
(except in \S~\ref{subsect:scanlontheorem}).
Any finite extension of $K$ will be
taken in $\overline{K}$.

An abelian variety $B$ over $\overline{K}$
is called a \emph{constant abelian variety} if there exists an abelian variety
$\widetilde{B}$ over $k$ such that 
$B = \widetilde{B} \otimes_k \overline{K}$
as abelian varieties.

Let $A$ be an abelian variety over $\overline{K}$.
A pair $\left( \widetilde{A}^{\overline{K}/k}, \Tr_A \right)$
consisting
of an abelian variety $\widetilde{A}^{\overline{K}/k}$ over $k$ and a homomorphism
$\Tr_A : \widetilde{A}^{\overline{K}/k} \otimes_{k} \overline{K} \to A$
of abelian varieties over $\overline{K}$
is called a \emph{$\overline{K}/k$-trace} of $A$ if
for any abelian variety $\widetilde{B}$ over $k$
and a homomorphism $\phi : \widetilde{B} \otimes_{k} \overline{K}
\to A$,
there exists a unique homomorphism $\Tr (\phi) :
\widetilde{B} \to 
\widetilde{A}^{\overline{K}/k} $ such that
$\phi$ factors as
$\phi = \Tr_A \circ \left( \Tr (\phi) \otimes_k \overline{K} \right)$.
It is unique by the universality,
and it is also known to exist.
We call $\Tr_A$ the \emph{trace homomorphism} of $A$.
See \cite{lang1} for more details.

Let $M_{K}$ be the set of points of $\mathfrak{B}$
of codimension $1$.
For any $v \in M_{K}$, the local ring $\OO_{\mathfrak{B}, v}$ 
is a discrete valuation
ring with fractional field $K$.
The order function on $\OO_{\mathfrak{B}, v}$
extends to a unique order function
$\ord_{v} : K^{\times} \to \ZZ$.
Recall that we have a fixed ample line bundle $\mathcal{H}$
on $\mathfrak{B}$.
Then we have
a
non-archimedean 
value $| \cdot |_{v , \mathcal{H}}$ on $K$
normalized
in such a way that
\begin{align*} 
| x |_{v , \mathcal{H}} := 
e^{- \deg_{\mathcal{H}} ( \overline{v} )
 \ord_{v} (x)}
\end{align*}
for any $x \in K^{\times}$,
where 
$\deg_{\mathcal{H}} ( \overline{v} )$
denotes the degree with respect to $\mathcal{H}$
of the closure $\overline{v}$ of $v$
in $\mathfrak{B}$.
It is well known that the set 
$
\{ | \cdot |_{v , \mathcal{H}} \}_{v \in M_{K}}$ 
of values satisfies the product formula,
and hence 
the notion of (absolute logarithmic) heights
with respect to this set of absolute values
is defined (cf. \cite[Chapter~3 \S~3]{lang2}).

We recall the notion of places of $\overline{K}$ and notation
introduced in \cite[\S~6.1]{yamaki6}.
For a finite extension $K'$ of $K$ in $\overline{K}$,
let $\mathfrak{B}'$ be the normalization of $\mathfrak{B}$ in $K'$
and
let $M_{K'}$ be the set of points of $\mathfrak{B}'$ of codimension $1$.
An element in $M_{K'}$ is called a \emph{place of $K'$}.
For a finite extension $K'' / K'$, we have a natural surjective map $M_{K''} \to M_{K'}$,
and thus we obtain a inverse system $(M_{K'})_{K'}$,
where $K'$ runs through the finite extensions of $K$ in $\overline{K}$.
Set $M_{\overline{K}} := \varprojlim_{K'} M_{K'}$.
We call an element of $M_{\overline{K}}$ a \emph{place of $\overline{K}$}.
Each $v \in M_{\overline{K}}$ gives a unique
absolute value on
$\overline{K}$ which extends $| \cdot |_{v_{K} , \mathcal{H}}$,
where $v_K$ is the image of $v$ by the canonical map $M_{\overline{K}}
\to M_{K}$.
We denote by $\overline{K}_{v}$ the completion 
of $\overline{K}$ with respect to that absolute value,
and we let
$\overline{K}_v^{\circ}$ denote the ring of integers of $\overline{K}_v$.

Let $A$ be an abelian variety over $\overline{K}$.
Let $v \in M_{\overline{K}}$ be a place.
We say that $A$ is \emph{non-degenerate at $v$}
if 
there exists an abelian scheme over $\overline{K}_v^{\circ}$
whose generic fiber equals $A$.
We say $A$ is \emph{nowhere degenerate}
if $A$ is non-degenerate at any $v \in M_{\overline{K}}$.
Those definitions are compatible with the terminology used in
\cite{yamaki6, yamaki7}.

We give a remark on our terminology of nowhere-degeneracy.
Suppose that $K'$ is a finite extension of $K$
and that $A'$ is an abelian variety over $K'$ with $A = A' \otimes_{K'}
\overline{K}$.
Let $v \in M_{\overline{K}}$ and let $v_{K'}$ be the image of $v$
by the natural map $M_{\overline{K}} \to M_{K'}$.
Then $A$ is non-degenerate at $v$ if and only if $A'$ has potentially 
good reduction at $v_{K'}$,
and thus $A$ is nowhere degenerate if and only if 
$A'$ has potentially 
good reduction everywhere.

Let $A$ be an abelian variety over $\overline{K}$
and let $L$ be an even ample line bundle.
Then there exists a unique height function $\widehat{h}_L$
on $A$ associated to $L$
such that $\widehat{h}_L$ is a quadratic form on the additive group 
$A \left( \overline{K} \right)$.
This is called the \emph{canonical height} associated to $L$.

\subsection{Canonical heights}

On an abelian variety,
we have a notion of canonical heights
not only for points but also
for positive dimensional cycles.
Let $L_0 , \ldots , L_d$ be line bundles 
on 
an abelian variety $A$
and
let $Z$ be a cycle of dimension $d$ on $A$.
Then 
we consider a real number called the
\emph{canonical height} of $Z$ with respect 
$L_0 , \ldots , L_d$,
which is denoted by 
$\widehat{h}_{L_0 , \ldots , L_d} (Z)$.
It is known that the assignment
$
(L_0 , \ldots , L_d , Z) \mapsto 
\widehat{h}_{L_0 , \ldots , L_d} (Z)
$
is multilinear.
When $L_0 = \cdots = L_d = L$ and
there is no danger of confusion,
we simply write
$\widehat{h}_{L} ( Z )$ instead of $\widehat{h}_{L_0 , \ldots , L_d} (Z)$.
We consider the canonical height of a closed subvariety $X$
of $A$ of pure dimension $d$ by regarding $X$ as
a cycle in a natural way.
We refer to \cite{gubler0, gubler2, gubler3} for more details.

The canonical heights of subvarieties are important
in the study of
density of small points because of the following proposition.

\begin{citeProposition} [Corollary~4.4 in \cite{gubler2}] \label{prop:dense-height0}
Let $A$ be an abelian variety over $\overline{K}$,
$L$ an even ample line bundle on $A$,
and
let $X$ be a
closed subvariety
of $A$.
Then
$X$ 
has dense small points if and only if
$\widehat{h}_{L} (X) = 0$.
\end{citeProposition}

In the case where
$A$ is a nowhere degenerate abelian variety,
the canonical height of a cycle can be
described 
in terms of intersection
products on models.
We refer to \cite{fulton} for intersection theory.
First, we recall the notion of models.
Let $X$ be a projective scheme over $\overline{K}$
and let $L$ be a line bundle on $X$.
Let $K'$ be a finite extension of $K$
and let $\mathfrak{B}'$ be  the normalization of $\mathfrak{B}$ in $K'$.
Let $\mathfrak{U}$ be an open subset of $\mathfrak{B}'$.
A proper morphism $f : \mathscr{X} \to \mathfrak{U}$
with geometric generic fiber $X$
is called a \emph{model} of $X$ over $\mathfrak{U}$.
Furthermore,
let $\mathscr{L}$ be a line bundle on $\mathscr{X}$
whose restriction to the
geometric generic fiber $X$ equals $L$.
Then the pair $(f, \mathscr{L})$ is called a \emph{model} of $(X,L)$
over $\mathfrak{U}$.

For a nowhere degenerate abelian variety $A$ over $\overline{K}$
and an even line bundle $L$ on $A$,
the following proposition
gives us a ``normalized'' model of $(A,L)$.
This is
a starting point to describe the canonical
height by intersection on models.

\begin{citeProposition} [Proposition~2.5 of \cite{yamaki7}]
\label{prop:model1}
Let $A$ be a nowhere degenerate abelian variety over $\overline{K}$
and let $L$ be a line bundle on $A$.
Then there exist
a finite extension $K'$ of $K$
that satisfies the following condition:
Let $\mathfrak{B}'$ be the normalization of $\mathfrak{B}$
in $K'$;
then there exist
an open subset $\mathfrak{U}$ of $\mathfrak{B}'$
with
$\codim ( \mathfrak{B}' 
\setminus \mathfrak{U} , \mathfrak{B}')
\geq 2$,
an abelian scheme $f : \mathscr{A} \to \mathfrak{U}$
with zero-section $0_f$,
and 
a line bundle $\mathscr{L}$ on $\mathscr{A}$
such that
$\left( f , \mathscr{L} \right)$ is a model of
$(A ,L)$
over $\mathfrak{U}$
with $0_{f}^{\ast} ( \mathscr{L} ) \cong \OO_{\mathfrak{U}}$.
\end{citeProposition}

Remark that $\mathfrak{U}$
as well as
$\mathfrak{B}'$
is normal
and that we have a natural finite 
surjective morphism $\mathfrak{B}' \to \mathfrak{B}$.

\begin{Remark} \label{rem:embedding}
We constructed in
\cite[Proposition~2.5]{yamaki7}
a model 
$\left( \bar{f} : \bar{\mathscr{A}} \to
\mathfrak{B}' , \bar{\mathscr{L}} \right)$
over $\mathfrak{B}'$ whose restriction over $\mathfrak{U}$ 
coincides with $(f , \mathscr{L})$
as in Proposition~\ref{prop:model1}.
Such an $\left( \bar{f} , \bar{\mathscr{L}} \right)$
is constructed from $(f , \mathscr{L})$
by using Nagata's embedding theorem.
We refer to \cite[Theorem~5.7]{vojta}
for a scheme-theoretic proof of this embedding theorem.
\end{Remark}

\begin{Remark}
\label{remark:new:nodeg-model}
Proposition~\ref{prop:model1}
says in particular that
if $A$ is nowhere degenerate, then 
there exist a finite extension $K'$ of $K$,
an open subset $\mathfrak{U}$ of $\mathfrak{B}'$
with $\codim ( \mathfrak{B}' \setminus \mathfrak{U} , \mathfrak{B}')
\geq 2$
where
$\mathfrak{B}'$ is the 
normalization
of $\mathfrak{B}$ in $K'$,
an abelian scheme $f : \mathscr{A} \to \mathfrak{U}$
with geometric generic fiber
$A$.
Remark that the converse of this also holds;
if $K'$ is a finite extension of $K$,
$\mathfrak{U}$ 
is an open subset of 
$\mathfrak{B}'$
with $\codim ( \mathfrak{B}' \setminus \mathfrak{U} , \mathfrak{B}')
\geq 2$
where $\mathfrak{B}'$ is the normalization of $\mathfrak{B}$ in $K'$,
and if $f : \mathscr{A} \to \mathfrak{U}$
is an abelian scheme,
then the geometric generic fiber $A$ of $f$ is a nowhere degenerate
abelian variety.
Indeed,
in this setting, 
since $\mathfrak{B}'$
is the normalization of $\mathfrak{B}$ in $K'$,
the set of places $M_{K'}$
equals the set of codimension $1$ points in $\mathfrak{B}'$.
Since
$\codim ( \mathfrak{B}' \setminus \mathfrak{U} , \mathfrak{B}')
\geq 2$,
that equals the set of codimension $1$ points in $\mathfrak{U}$.
It follows that the generic fiber 
of $f$ is an abelian variety
that has good reduction at any $v_{K'} \in M_{K'}$,
and thus
$A$ is nowhere degenerate.
\end{Remark}

Using a model
as in Proposition~\ref{prop:model1},
we describe the canonical height of a cycle
in terms of intersection
by
\cite[Theorem~3.5~(d)]{gubler3}.
We show the following lemma, which will 
follow easily from that theorem,
where we recall $b := \dim ( \mathfrak{B})$
and $\mathcal{H}$ is a fixed ample line bundle on $\mathfrak{B}$.

\begin{Lemma} \label{lem:height-intersection1}
Let $A$ be a nowhere degenerate abelian variety over $\overline{K}$
and let $L$ be an even line bundle on $A$.
Let $X$ be a closed subscheme of $A$
of pure dimension $d$.
Let 
$K'$,
$\mathfrak{B}'$,
$\mathfrak{U}$,
$f ; \mathscr{A} \to \mathfrak{U}$,
$0_f$, and $\mathscr{L}$
be as in Proposition~\ref{prop:model1}.
Let $\mathcal{H}'$ be the pull-back of $\mathcal{H}$ by the finite morphism
$\mathfrak{B}' \to \mathfrak{B}$.
Let $\mathscr{X}$ be the closure of $X$ in $\mathscr{A}$.
Let $\mathfrak{D}$ be a $(b-1)$-cycle on $\mathfrak{B}'$
such that 
$[\mathfrak{D} \cap \mathfrak{U}]
= f_{\ast}
\left(
\cherncl_1 ( \mathscr{L})^{\cdot (d+1)}
\cdot
[\mathscr{X}] 
\right)$
as cycle classes
on $\mathfrak{U}$.
Then
\[
\widehat{h}_L (X)
= 
\frac{
\deg_{\mathcal{H}'} 
[\mathfrak{D}] }{ [K' : K]}
.
\]
\end{Lemma}

\Proof
As in Remark~\ref{rem:embedding},
we take 
a proper morphism 
$\bar{f} : \bar{\mathscr{A}} \to \mathfrak{B}'$
and a line bundle $\bar{\mathscr{L}}$
such that $f$ is the pull-back of $\bar{f}$ by 
the open immersion
$\mathfrak{U} \hookrightarrow \mathfrak{B}'$
and $\bar{\mathscr{L}}|_{\mathscr{A}} = \mathscr{L}$.
Let $\bar{\mathscr{X}}$ be the closure of $\mathscr{X}$
in
$\bar{\mathscr{A}}$.

Note that
\addtocounter{Claim}{1}
\begin{align} \label{eq:compactheight}
\widehat{h}_L ( X )
=
\frac{
\deg_{\mathcal{H}'}
\bar{f}_{\ast}
\left( \cherncl_1 \left( \bar{\mathscr{L}} \right)^{\cdot (d+1)} \cdot
\left[ 
\bar{\mathscr{X}}
\right]
\right)}
{[K' : K]}
.
\end{align}
Indeed,
since $\mathfrak{U}$ is normal,
$f : \mathscr{A} \to \mathfrak{U}$
is projective by \cite[Th\'eor\`eme~XI~1.4]{raynaud0}.
By \cite[Corollaire~5.7.14]{RG}
or by \cite[Corollary~2.6]{vojta},
there exists a proper morphism $\mu : \bar{\mathscr{A}}^{\dagger} \to \bar{\mathscr{A}}$
isomorphic over $\mathscr{A}$
such that $\bar{f}^{\dagger} := \bar{f} \circ \mu$
is projective and such that $\mu^{-1} ( \mathscr{A} )$ is scheme-theoretically
dense in $\bar{\mathscr{A}}^{\dagger}$.
Noting that
 $\bar{f}^{\dagger} : \bar{\mathscr{A}}^{\dagger} \to \mathfrak{B}'$
is also a model of $A$,
we let $\bar{\mathscr{X}}^{\dagger}$ be the closure of 
$X$ in $\bar{\mathscr{A}}^{\dagger}$.
Since $\bar{f}^{\dagger}$ is projective,
\cite[Theorem~3.5~(d)]{gubler3}
gives us
\begin{align*} \label{eq:compactheight}
\widehat{h}_L ( X )
=
\frac{
\deg_{\mathcal{H}'}
\bar{f}^{\dagger}_{\ast}
\left( \cherncl_1 \left( \mu^{\ast} \left(
\bar{\mathscr{L}} \right) \right)^{\cdot (d+1)} \cdot
\left[ 
\bar{\mathscr{X}}^{\dagger}
\right]
\right)}
{[K' : K]}
.
\end{align*}
By the projection formula,
we thus obtain (\ref{eq:compactheight}).

Since $\bar{f}_{\ast}
\left( \cherncl_1 \left( \bar{\mathscr{L}} \right)^{\cdot (d+1)} \cdot
\left[ 
\bar{\mathscr{X}}
\right]
\right)
$
is
a cycle class on $\mathfrak{B}'$ of
codimension $1$
whose restriction to $\mathfrak{U}$
equals
$
f_{\ast}
\left(
\cherncl_1 ( \mathscr{L})^{\cdot (d+1)}
\cdot
[\mathscr{X}] 
\right)$
and 
since
$\codim \left( \mathfrak{B}' \setminus \mathfrak{U} , \mathfrak{B}'
\right) \geq 2$,
we have 
\[
[\mathfrak{D}] =
\bar{f}_{\ast}
\left( \cherncl_1 \left( \bar{\mathscr{L}} \right)^{\cdot (d+1)} \cdot
\left[ 
\bar{\mathscr{X}}
\right]
\right)
\]
as cycle classes.
Thus the lemma follows from (\ref{eq:compactheight}).
\QED

\begin{Remark} \label{rem:heightabelianvariety}
Let $\bar{f} : \bar{\mathscr{A}} \to \mathfrak{B}'$
be as in Remark~\ref{rem:embedding}.
Then 
\[
\deg_{\mathcal{H}'}
\bar{f}_{\ast}
\left( \cherncl_1 \left( \bar{\mathscr{L}} \right)^{\cdot (\dim (A)+1)} \cdot
\left[ 
\bar{\mathscr{A}}
\right]
\right) = 0.
\]
Indeed,
since $A$ has dense small points, it has canonical height $0$
(cf. Proposition~\ref{prop:dense-height0}),
and hence the equality follows from Lemma~\ref{lem:height-intersection1}
(or (\ref{eq:compactheight})).
\end{Remark}

\section{Bertini-type theorem} \label{sect:Bertini}

The purpose of this section is to show
Proposition~\ref{prop:bertini5}, which 
concerns 
curves 
which are intersection of general hyperplanes
on a projective variety.
This proposition will be applied later to a finite covering
of $\mathfrak{B}$
in the case of $b := \dim ( \mathfrak{B} ) \geq 2$.

Let $\OO (1)$ denote the tautological line bundle of
$\PP^{N}_k$
and let
$|\OO (1)|$ denote
the complete linear system of hyperplanes
on $\PP^{N}_k$.
Let $X$ be a subvariety of $\PP^{N}_k$.
Let $r$ be a natural number.
For any $\underline{H} = \left(
H_1 , \ldots , H_{r} 
\right)
\in |\OO (1)|^{r} 
$,
we write
$X_{\underline{H}} := X \cap H_1 \cap  \cdots \cap H_{r}$.

\begin{Lemma} \label{lem:bertini3}
Let $X$ be an irreducible projective scheme 
of dimension $d$
over $k$
with a 
closed embedding $X \subset \PP^{N}_k$.
Suppose that there exists 
a regular open subscheme
$U$ with
$\codim ( X \setminus U , X) 
\geq 2$.
Let $r$ be a positive integer with $r \leq d - 1$.
Then
there exists 
a dense open subset 
$V (k)
 \subset 
|\OO (1)|^{r} $
such that any 
$\underline{H} = \left(
H_1 , \ldots , H_{r} 
\right)
\in V (k)$
satisfies
the following:
\begin{enumerate}
\renewcommand{\labelenumi}{(\alph{enumi})}
\item
$X_{\underline{H}}$ is irreducible
and has dimension $d-r$;
\item
$U_{\underline{H}}$ is regular;
\item
$\codim (X_{\underline{H}} \setminus U , X_{\underline{H}}) \geq 2$.
\end{enumerate}
\end{Lemma}

\Proof
Note that the assignment
$
 \left(
H_1 , \ldots , H_{r} 
\right) \mapsto 
H_1 \cap \cdots \cap H_{r} 
$
defines a dominant rational map from 
$|\OO (1)|^{r} $
to the grassmannian variety 
of codimension $r$ linear subspaces of $\PP^{N}_k$.
Then the lemma follows from \cite[Corollary~6.11]{jouanolou}
immediately.
Indeed,
by \cite[Corollary~6.11]{jouanolou},
there exists a dense open subset $V_1
(k)$ of 
$|\OO (1)|^{r} $
such that for any $\underline{H} \in V_1 (k)$,
$X_{\underline{H}}$ is irreducible
subscheme of dimension $d-r$
and
$U_{\underline{H}}$ is regular.
Further by \cite[Corollary~6.11~(1)~(b)]{jouanolou},
there exists a
dense open subset $V_2 (k)$ of 
$|\OO (1)|^{r} $
such that for any $\underline{H} \in V_1 (k)$,
$\codim \left( (X \setminus U)_{\underline{H}} , X \right) \geq 2+r$.
Putting $V (k) := V_1 (k) \cap V_2 (k)$,
we see that any $\underline{H}
\in V (k)$ satisfies conditions (a), (b), and (c).
\QED

We sometimes have to consider the linear system
consisting of hyperplanes which pass through a
specified point.
For any $x \in X (k)$,
we set
\[
|\OO (1)|^{r}_x
:=
\{
(H_1 , \ldots , H_r )
\in |\OO (1)|^{r}_x
\mid
x \in H_1 \cap \cdots \cap H_r
\}
.
\]
Further,
for an open subset $V (k)$ of $|\OO (1)|^{r}$,
we set
$V_{x}(k) := 
V (k) \cap | \OO (1) |^{r}_x$.

\begin{Lemma} \label{lem:generalmember2-1}
Let $d$ be an integer with $d \geq 2$.
Let $r$ be a positive integer with $r \leq d - 1$
and let $V (k)$ be a dense open subset of 
$|\OO (1)|^{r} $.
Let $x \in X(k)$ be a point
with
$
V_{x}(k) 
\neq \emptyset
$.
Suppose that $W$  is a
non-empty closed subset of 
$\PP^{N}_k$ with $\dim (W) \leq d-1$.
Then
there exists an
$\underline{H} = (H_1 , \ldots , H_r )\in V_{x} (k)$ such that
$\dim \left( W_{\underline{H}} \right) \leq d-r-1$.
(Remark that the dimension of a closed subset is
the maximum of the dimensions of its irreducible components.
Remark also that $\dim ( \emptyset ) := -1$ by convention.)
\end{Lemma}

\Proof
We prove the lemma by induction on $r$.
First let $r=1$.
If $\dim (W) \leq d - 2$,
then $\dim (W \cap H) \leq d-2$
for any hyperplane $H$.
Therefore we may assume $\dim (W) = d -1$.
Let $q_{1} , \ldots , q_m$ be the generic points
of the irreducible components of $W$ of dimension $d-1$.
Since $d-1 \geq 1$, no $q_i$ equals $x$.
Therefore
\[
U(k):=
\{ H \in |\OO (1)|_x 
\mid
q_1 \notin H , \ldots ,  q_m \notin H \}
\]
is a dense open subset of $|\OO (1)|_x$,
and hence there exists an 
$H \in V_x (k)  \cap U(k)$.
We then have $\dim (W \cap H) \leq d - 2$,
and thus
we have the lemma for $r = 1$.

Suppose that we have the assertion up to $r$ ($1 \leq r \leq d-2$),
and we are going to show it for $r+1$.
Let $p : |\OO (1)|_x^{r+1} \to |\OO(1)|_x^{r}$ the map
given by $p (H_1 , \ldots , H_r , H_{r+1} ) = (H_1 , \ldots , H_r  )$.
Then the image
$p ( V_x(k) ) $ is a dense open subset of $|\OO(1)|_x^{r}$.
By the induction hypothesis,
there exists an $\underline{H} = (H_1 , \ldots , H_r) \in p ( V_x(k) ) $ 
such that
$\dim \left( W_{\underline{H}} \right) \leq d - r - 1$.
Since $p^{-1} ( \underline{H} ) \cap V_x(k) \neq \emptyset$,
this set is a dense open subset of 
$\{ \underline{H} \} \times | \OO (1)|_x$.
Applying the lemma for $r=1$ to
$W_{\underline{H}}$,
we obtain an $H_{r+1} \in p^{-1} ( \underline{H} ) \cap V_x(k)$
such that $\dim \left( W_{\underline{H}} \cap H_{r+1}
\right) \leq d - r -2$.
This shows the assertion for $r+1$.
Thus we obtain the lemma.
\QED

Using Lemma~\ref{lem:generalmember2-1},
we obtain the following:

\begin{Lemma} \label{lem:generalmember2}
Let $X$ be an irreducible projective scheme
of dimension $d$
over $k$
with a 
closed embedding $X \subset \PP^{N}_k$.
Let $r$ be a positive integer with $r \leq d - 1$
and let $V (k)$ be a dense open subset of 
$|\OO (1)|^{r} $.
Further, let $x \in X(k)$ be a point
such that
$
V_{x}(k) \neq \emptyset
$.
Then for any proper closed subset $W$ of $X$,
there exists an
$\underline{H} = (H_1 , \ldots , H_r )\in V_{x} (k)$ such that
$\Supp ( 
X_{\underline{H}} ) \nsubseteq W$.
\end{Lemma}

\Proof
Remark that $d \geq r+ 1 \geq 2$.
Remark also that $\dim \left( X_{\underline{H}} \right) \geq d - r$ for any
$\underline{H} \in |\OO (1)|^{r}$.
Let $W$ be any proper closed subset of $X$.
Then $\dim (W) \leq d- 1$.
By Lemma~\ref{lem:generalmember2-1},
there exists an 
$\underline{H} = (H_1 , \ldots , H_r )\in V_{x} (k)$ such that
$\dim \left( W_{\underline{H}} \right) \leq d - r - 1$.
Since $\dim \left( X_{\underline{H}} \right) \geq d - r$,
we have $X_{\underline{H}}
\nsubseteq W_{\underline{H}}$,
and thus
$\Supp ( 
X_{\underline{H}} ) \nsubseteq W$.
\QED

We show one more lemma.

\begin{Lemma} \label{lem:cl:prop:bertini2}
Let $X$ be an irreducible closed subscheme 
of dimension $d$
of $
\PP^{N}_k$.
Suppose that there exists a regular open subscheme
$U$
of 
$X$ 
with
$\codim ( X \setminus U , X) 
\geq 2$.
Let $r$ be a positive integer with $r \leq d - 1$
and let
$V
(k) 
$ be a dense open subset
of
$|\OO (1)|^{r} $.
Then 
there exists a dense open subset $U' \subset X$
such that for any $x \in U' (k)$,
there exists an $( H_1 , \ldots , H_r ) \in V(k)$
with $x \in H_1 \cap \cdots \cap H_r$.
\end{Lemma}

\Proof
Set
\[
\mathcal{Z} (k) 
:=
\{ ( x , H_1 , \ldots , H_r ) \in X (k) \times V(k)
\mid 
x \in H_1 \cap \cdots \cap H_r \}
,
\]
which is a non-empty closed subset 
of $X (k) \times V(k)$.
Let 
$g : \mathcal{Z} (k) \to X (k)$ be the
restriction of the canonical projection 
$X (k)  \times V (k) \to X (k)$.
We take an $x \in g ( \mathcal{Z} (k))$.
Note that 
\[
V_{x} (k)
:=
\{
\underline{H} = (H_1 , \ldots , H_r)
\in V (k)
\mid
x \in H_1 \cap \cdots \cap H_r
\}
\neq \emptyset
.
\]

Let $W (k)$ be the closure of $g ( \mathcal{Z} (k))$ in $X (k)$.
Then
we claim $W(k) = X (k)$
by contradiction.
Indeed,
if this is not the case,
then $W(k) \subsetneq X (k)$,
and hence
Lemma~\ref{lem:generalmember2} gives us 
an $\underline{H} = (H_1 , \ldots , H_r) \in V_x(k)$
such that $X_{\underline{H} } (k) 
\nsubseteq W (k)$.
Now we can take
a point $x'$ of $X_{\underline{H}} (k)
\setminus W (k)$.
Then
$x' \in g ( \mathcal{Z} (k) )$ by the definition of $\mathcal{Z} (k)$.
It follows that $x' \in g ( \mathcal{Z} (k) ) \setminus W (k)
\subset W (k) \setminus W (k) = \emptyset$,
which is a contradiction.

Therefore $g : \mathcal{Z}(k) \to X (k)$
is a dominant morphism.
By Chevalley's theorem, 
there exists a non-empty open subset $U'$ of $X$ such that 
$U'(k) \subset g ( \mathcal{Z}(k) )$.
Then $U'$ satisfies the required condition,
and this completes the proof.
\QED

As a consequence of the arguments so far,
we have the following proposition.

\begin{Proposition} \label{prop:bertini5}
Let $X$ be an irreducible projective scheme
of dimension $d$
over $k$
and let
$\mathscr{L}$
be
a very ample line bundle on $X$.
Let $|\mathscr{L}|$ denote the complete linear system associated
to $\mathscr{L}$.
Suppose that there exists a 
regular open subscheme
$U$ of $X$
such that
$\codim ( X \setminus U , X) 
\geq 2$.
Then the following hold.
\begin{enumerate}
\item
There exists 
a dense open subset $V (k) \subset |\mathscr{L}|^{d - 1}$
such that for any 
$
\left(
D_1 , \ldots , D_{d - 1} 
\right)
\in V (k)$,
$C:= D_1 \cap \cdots \cap D_{d - 1}$
is an irreducible projective non-singular curve contained in $U$.
\item
Let $V(k)$ be as in (1).
For any $x \in X (k)$,
set
\[
|\mathscr{L}|^{d-1}_x
:=
\left\{ 
\left.
\left(
D_1 , \ldots , D_{d - 1} 
\right) \in | \mathscr{L} |^{d-1}
\ 
\right|
x \in D_1 \cap \cdots \cap D_{d - 1}
\right\}
\]
and $V_{x} (k)
:=
V \cap |\mathscr{L}|^{d-1}_x$.
Then there exists a dense open subset $U' \subset U$
such that for any $x \in U' (k)$, 
$V_{x} (k)$
is a dense open subset of $|\mathscr{L}|^{d-1}_x$.
\item
Furthermore, let $U'$ be as in (2).
Take any $x \in U' (k)$.
Then for any proper closed subset $W$ of $ X$,
there exists a $\left(
D_1 , \ldots , D_{d - 1} 
\right)
\in V_x (k)$
such that
$\Supp (D_1 \cap \cdots \cap D_{d - 1}) \nsubseteq W$.
\end{enumerate}
\end{Proposition}

\Proof
Let $X \hookrightarrow \PP^{N}_k$ be the closed embedding
associated to the complete linear system $| \mathscr{L} |$.
Then we have $| \mathscr{L} | = | \OO (1)|$,
where $\OO (1)$ is the tautological line bundle of $\PP^{N}_k$.
Now, we see that 
(1) follows from Lemma~\ref{lem:bertini3}
for $r = d-1$.
Since $V_x (k)$ is open in $|\mathscr{L}|^{d-1}_x$,
(2) follows from Lemma~\ref{lem:cl:prop:bertini2}.
Finally, (3) follows from Lemma~\ref{lem:generalmember2}.
\QED

\section{Positivity of the canonical height} \label{sect:positivity}

In this section, we prove Proposition~\ref{prop:main1},
which claims the positivity of the canonical height 
of an effective ample divisor on a
nowhere degenerate abelian variety
over $\overline{K}$ with trivial 
$\overline{K}/k$-trace.

Let $\mathscr{A} \to S$ is an abelian scheme.
We use the following notation.
For a section $\sigma$ of it,
we have the translate morphism 
by $\sigma$,
which is denoted by
$T_{\sigma} : \mathscr{A} \to \mathscr{A}$.
For an $n$-times endomorphism $[n] : \mathscr{A} \to \mathscr{A}$,
we set $\mathscr{A} [n] := \Ker [n]$.
That is a finite scheme over $S$.

\subsection{Translates of effective divisors on abelian varieties}

In this subsection, 
let $A$ be an abelian variety over $k$
and let $D$ be an effective divisor on $A$.
We show two lemmas concerning translates of $D$.

\begin{Lemma} \label{lem:notinsupport}
For any $a \in A (k)$,
there exists a 
dense open subset $V_a \subset A$ such that 
$a \notin \Supp  \left( T_{\tau}^{\ast}( D )
+
T_{-\tau}^{\ast}( D ) \right)$
for any
$\tau \in V_a ( k )$,
where ``$+$'' means the sum of the divisors.
\end{Lemma}

\Proof
We see that $a \in \Supp  \left( T_{\tau}^{\ast}( D ) \right)$
is equivalent to $a \in D - \tau$,
and this is equivalent to $\tau \in \Supp  \left( T_{a}^{\ast}( D )
\right)$.
This is a closed condition for $\tau \in A (k)$.
Similarly, 
we see that $a \in \Supp  \left( T_{-\tau}^{\ast}( D )
\right)$
is a closed condition for $\tau$.
Hence
$a \in \Supp  \left( T_{\tau}^{\ast}( D )
\right)
\cup
\Supp  \left( T_{-\tau}^{\ast}( D )
\right)$
is an open condition for $\tau$.
Since $\Supp  \left( T_{a}^{\ast}( D ) + 
T_{-\tau}^{\ast}( D )
\right) 
\subset
\Supp  \left( T_{\tau}^{\ast}( D )
\right)
\cup
\Supp  \left( T_{-\tau}^{\ast}( D )
\right)
\subsetneq A$,
that shows the existence of
a desired dense open subset $V_a$.
\QED

\begin{Lemma} \label{lem:emptyintersctiononabelianvariety}
Let $l$ be a prime number with $l \neq \ch (k)$.
Then there exists an $m \in \NN$
such that
\[
\bigcap_{ \tau \in A [ l^{m} ]} 
\Supp  \left( T_{\tau}^{\ast}( D ) + T_{-\tau}^{\ast} ( D ) \right)
=
\emptyset
.
\]
\end{Lemma}

\Proof
For each $n \in \NN$,
we put
\[
S_n
:=
\bigcap_{ \tau \in A [ l^{n} ]} 
\Supp  \left( T_{\tau}^{\ast}( D ) + T_{-\tau}^{\ast} ( D ) \right)
.
\]
Since $(S_n)_{ \in \NN}$ is a descending sequence of closed subsets of $A$,
there exists
an $m \in \NN$
such that $S_n = S_m$ for any $n \geq m$.
What we should show is that $S_m = \emptyset$.
To show this by contradiction,
suppose that $S_m \neq \emptyset$.
Since $S_m$ is closed,
we then take a closed point
$s \in S_m (k)$.
Let $V_{s}$ be a dense open subset of $A$ as in Lemma~\ref{lem:notinsupport}.
Since $l \neq \ch (k)$,
$\bigcup_{n \in \NN} A [l^{n}]$
is dense in $A$, 
and hence we take
a point
$
\tau_{0}
\in
V_{s} (k)
\cap 
\left(
\bigcup_{n \in \NN} A [l^{n}]
\right) 
$.
Then 
$s \notin 
\Supp  \left( T_{\tau_0}^{\ast}( D ) + T_{-\tau_0}^{\ast} ( D ) \right)$.
On the other hand,
there exists
an 
$m_0 \in \NN$ 
with $m_0 \geq m$ such that $\tau_0 \in A [l^{m_0}]$.
Since $s \notin 
\Supp  \left( T_{\tau_0}^{\ast}( D ) + T_{-\tau_0}^{\ast} ( D ) \right)$,
it follows that
\[
s
\notin 
\bigcap_{ \tau \in A [ l^{m_0} ]} 
\Supp  \left( T_{\tau}^{\ast}( D ) + T_{-\tau}^{\ast} ( D ) \right)
=
S_{m_0}.
\]
However, since we have taken $m$ and $s$ so that $S_{m_0}
= S_m$ and $s \in S_m$,
that
is a contradiction.
Thus we conclude that $S_{m} = \emptyset$,
which completes the proof of the lemma.
\QED

\subsection{Base-point freeness on abelian schemes}

The purpose of this subsection 
is to establish Proposition~\ref{prop:basepointfree},
the base-point freeness 
on an abelian scheme
of an effective even 
line bundle
which is normalized to be trivial along the zero-section.

First, we show the following lemma,
which generalizes 
Lemma~\ref{lem:emptyintersctiononabelianvariety}
on abelian schemes.

\begin{Lemma} \label{lem:bpf1}
Let $\mathfrak{U}$ be an
irreducible
noetherian scheme,
$f : \mathscr{A} \to \mathfrak{U}$ an
abelian scheme,
and let
$\mathscr{D}$ be an 
effective Cartier divisor on $\mathscr{A}$
flat over $\mathfrak{U}$.
Let $l$ be a prime number.
Suppose that 
$l$ does not equal to the characteristic of
the residue field at
any $u \in \mathfrak{U}$.
Then
there exist an $m \in \NN$ and
 a finite \'etale morphism 
$\mathfrak{U}' \to \mathfrak{U}$ 
with $\mathfrak{U}'$ irreducible
that satisfy the following conditions:
Let $f' : \mathscr{A}' \to \mathfrak{U}'$
and
$\mathscr{D}'$
be the base-change of 
$f$ and $\mathscr{D}$,
respectively
by the morphism $\mathfrak{U}' \to \mathfrak{U}$;
then
we have
\addtocounter{Claim}{1}
\begin{align} \label{eq:condition-lemmabpf1}
\bigcap_{\tau \in \mathscr{A}' [ l^{m} ] ( \mathfrak{U}' )}
\Supp
\left( T_{\tau }^{\ast} ( \mathscr{D}' )
+
T_{- \tau}^{\ast} ( \mathscr{D}')
\right)
=
\emptyset
,
\end{align}
where 
$\mathscr{A}' [ l^{m} ] ( \mathfrak{U}' )$ is the group of sections of 
$\mathscr{A}' [ l^{m} ] \to \mathfrak{U}'$.
\end{Lemma}

\Proof
We first construct 
by induction on $n \in \NN$ a sequence 
$\{ q_n : \mathfrak{U}_n \to \mathfrak{U} \}_{n \in \NN}$
of finite \'etale morphisms with $\mathfrak{U}_n$
irreducible 
such that
$\mathscr{A}_n
[ l^{n} ] 
=
\coprod_{\tau \in \mathscr{A}_n [ l^{n} ] (\mathfrak{U}_n)}
\tau
\left(
\mathfrak{U}_n
\right)$
where 
$f_n : \mathscr{A}_n
\to \mathfrak{U}_n
$ is the base-change of $f$ by $q_n$.
Since 
there does not exist a residue field 
with characteristic
$l$
on $\mathfrak{U}$,
the morphism
$\mathscr{A} [ l^{n} ] \to \mathfrak{U}$
is finite and \'etale for any $n \in \NN$.
It follows that there exists
a finite \'etale
morphism $q_1 : \mathfrak{U}_{1} \to \mathfrak{U}$
with $\mathfrak{U}_{1}$ irreducible
such that 
$\mathscr{A}_1
[ l ] 
=
\coprod_{\tau \in \mathscr{A}_1 [ l ] (\mathfrak{U}_1)}
\tau
\left(
\mathfrak{U}_1
\right)$.
Let $m$ be a natural number and
suppose that we have constructed a sequence
$\{ q_n : \mathfrak{U}_n \to \mathfrak{U} \}_{n = 1 ,\ldots , m}$
satisfying the condition.
Then by the same argument as above,
there exists a finite \'etale morphism $r : \mathfrak{U}_{m+1} \to
\mathfrak{U}_m$ 
with $\mathfrak{U}_{m+1}$ irreducible
such that
$\mathscr{A}_{m+1}
[ l^{m+1} ] 
=
\coprod_{\tau \in \mathscr{A}_{m+1} [ l^{m+1} ] (\mathfrak{U}_{m+1})}
\tau
\left(
\mathfrak{U}_{m+1}
\right)$.
Let $q_{m+1} : \mathfrak{U}_{m+1} \to \mathfrak{U}$ be the composite
$q_m \circ r$.
Then $q_{m+1}$ satisfies the required condition.
Thus we obtain $\{ q_n : \mathfrak{U}_n \to \mathfrak{U} \}_{n \in \NN}$ by
induction.

Let 
$\mathscr{D}_n$ 
be the pull-back of 
$\mathscr{D}$ 
by $\mathscr{A}_n
\to \mathscr{A}$
and
put
\[
\mathscr{F}_n :=
\bigcap_{\tau
\in \mathscr{A}_n [ l^{n} ] ( \mathfrak{U}_n )}
\Supp
\left( T_{\tau }^{\ast} ( \mathscr{D}_n )
+
T_{- \tau}^{\ast} ( \mathscr{D}_n)
\right)
.
\]
Then $\mathfrak{F}_n : =
q_{n} \left( f_n \left( \mathscr{F}_n \right) \right)$
is a closed subset of $\mathfrak{U}$.
Noting that there exists a natural map
$\mathscr{A}_n [ l^{n} ] ( \mathfrak{U}_n )
\hookrightarrow
\mathscr{A}_{n+1} [ l^{n+1} ] ( \mathfrak{U}_{n+1} )$,
we see that $\mathfrak{F}_n \supset \mathfrak{F}_{n + 1}$ for any $n \in \NN$.
Since $\mathfrak{U}$ is noetherian,
there exists an $m \in \NN$ such that 
for any $n \geq m$,
$\mathfrak{F}_n = \mathfrak{F}_{m}$
holds.

It then suffices to show
$\mathfrak{F}_{m} = \emptyset$
for this lemma;
Indeed, if this holds,
letting
$\mathfrak{U}' \to \mathfrak{U}$ be the morphism
$q_m : \mathfrak{U}_{m} \to \mathfrak{U}$,
we then see that
$m$ and this $\mathfrak{U}' \to \mathfrak{U}$
satisfy the required condition in the lemma.
We prove $\mathfrak{F}_{m} = \emptyset$
by contradiction;
Suppose that $\mathfrak{F}_{m} \neq \emptyset$.
Then there exists a geometric point $\overline{u} \in \mathfrak{F}_{m}
\subset \mathfrak{U}$.
Since $\mathscr{D}$ is
an effective Cartier divisor flat over $\mathfrak{U}$,
$\mathscr{D} \cap f^{-1} (\overline{u})$
is an effective divisor on 
the abelian variety
$f^{-1} (\overline{u})$.
By Lemma~\ref{lem:emptyintersctiononabelianvariety},
there exists an $m_0 \in \NN$ 
with $m_0 \geq m$ such that
\[
\bigcap_{\tau
\in f^{-1} (\overline{u}) \left[ l^{m_0} \right]}
\Supp
\left( T_{\tau }^{\ast} \left( \mathscr{D} \cap f^{-1} (\overline{u}) \right)
+
T_{- \tau}^{\ast} \left( \mathscr{D} \cap f^{-1} (\overline{u}) \right)
\right)
= \emptyset
.
\]
For any $\overline{u}_{m_0} \in q_{m_0}^{-1} ( \overline{u} )$,
remark that
the natural morphism
$\mathscr{A}_{m_0} \to \mathscr{A}$
restricts to isomorphism
$f_{m_0}^{-1} 
(\overline{u}_{m_0}) \cong
f^{-1} (\overline{u})$,
and via that isomorphism,
we have $\mathscr{D}_{m_0} \cap f_{m_0}^{-1} 
(\overline{u}_{m_0})
=
\mathscr{D} \cap f^{-1} ( \overline{u} )$.
Then
we have
\addtocounter{Claim}{1}
\begin{align} \label{eq:m_0}
\bigcap_{\tau
\in f_{m_0}^{-1} (\overline{u}_{m_0}) \left[ l^{m_0} \right]}
\Supp
\left( T_{\tau }^{\ast} \left( \mathscr{D}_{m_0} \cap f_{m_0}^{-1} 
(\overline{u}_{m_0}) \right)
+
T_{- \tau}^{\ast} \left( \mathscr{D}_{m_0} \cap f_{m_0}^{-1} (\overline{u}_{m_0}) 
\right)
\right)
= \emptyset
.
\end{align}
Since 
$\mathscr{A}_{m_0}
[ l^{m_0} ] 
=
\coprod_{\tau \in \mathscr{A}_{m_0} [ l^{m_0} ] (\mathfrak{U}_{m_0})}
\tau
\left(
\mathfrak{U}_{m_0}
\right)$
and the residue field of $\overline{u}_{m_0}$
does not equal $l$,
the natural map 
$\mathscr{A}_{m_0} [ l^{m_0} ]
\left(
\mathfrak{U}_{m_0}
\right) \to
f_{m_0}^{-1} 
(\overline{u}_{m_0}) \left[ l^{m_0} \right]$
given by restriction 
is an isomorphism.
It follows that
\begin{align*}
&
\bigcap_{\tau
\in f_{m_0}^{-1} (\overline{u}_{m_0}) \left[ l^{m_0} \right]}
\Supp
\left( T_{\tau }^{\ast} \left( \mathscr{D}_{m_0} \cap f_{m_0}^{-1} 
(\overline{u}_{m_0}) \right)
+
T_{- \tau}^{\ast} \left( \mathscr{D}_{m_0} \cap f_{m_0}^{-1} (\overline{u}_{m_0}) 
\right)
\right)
\\
&=
\bigcap_{\tilde{\tau}
\in \mathscr{A}_{m_0} [ l^{m_0} ] ( \mathfrak{U}_{m_0} )}
\Supp
\left( T_{\tilde{\tau} }^{\ast} ( \mathscr{D}_{m_0} )
+
T_{- \tilde{\tau}}^{\ast} ( \mathscr{D}_{m_0})
\right)
\cap
f_{m_0}^{-1} ( \overline{u}_{m_0} ).
\end{align*}
Therefore (\ref{eq:m_0})
shows that
\[
\bigcap_{\tilde{\tau}
\in \mathscr{A}_{m_0} [ l^{m_0} ] ( \mathfrak{U}_{m_0} )}
\Supp
\left( T_{\tilde{\tau} }^{\ast} ( \mathscr{D}_{m_0} )
+
T_{- \tilde{\tau}}^{\ast} ( \mathscr{D}_{m_0})
\right)
\cap
f_{m_0}^{-1} ( \overline{u}_{m_0} )
=
\emptyset
\]
for any
$\overline{u}_{m_0} \in q_{m_0}^{-1} ( \overline{u} )$.
On the other hand,
since $\overline{u} \in \mathfrak{F}_{m} = \mathfrak{F}_{m_0}$,
there exists an 
$\overline{u}_{m_0}' \in q_{m_0}^{-1} ( \overline{u} )$
such that
\[
\bigcap_{\tilde{\tau}
\in \mathscr{A}_{m_0} [ l^{m_0} ] ( \mathfrak{U}_{m_0} )}
\Supp
\left( T_{\tilde{\tau} }^{\ast} ( \mathscr{D}_{m_0} )
+
T_{- \tilde{\tau}}^{\ast} ( \mathscr{D}_{m_0})
\right)
\cap
f_{m_0}^{-1} ( \overline{u}_{m_0}' )
=
\mathscr{F}_{m_0} \cap 
f_{m_0}^{-1} ( \overline{u}_{m_0}' )
\neq
\emptyset
.
\]
However,
that is a contradiction.
Thus we conclude $\mathfrak{F}_{m} = \emptyset$,
which completes the proof of the lemma.
\QED

Now we show the following base-point freeness
on a model.

\begin{Proposition} \label{prop:basepointfree}
Let $K'$ be a finite extension of $K$,
$\mathfrak{B}'$ the normalization
of $\mathfrak{B}$ in $K'$,
$\mathfrak{U}$
an open subset of $\mathfrak{B}'$
with
$\codim \left( \mathfrak{B}' \setminus
\mathfrak{U} , \mathfrak{B}' \right) \geq 2$,
and
let $f : \mathscr{A} \to \mathfrak{U}$ be an
abelian scheme
with zero-section $0_f$.
Let $\mathscr{D}$ be an effective 
Cartier divisor on $\mathscr{A}$
that
is flat over $\mathfrak{U}$.
Suppose that
the restriction of
$\OO_{\mathscr{A}} ( \mathscr{D} )$ to the generic fiber 
of $f$ is an even line bundle and
that
$0_f^{\ast} \left( \OO_{\mathscr{A}} ( \mathscr{D} ) \right)
\cong \OO_{\mathfrak{U}}$.
Then
there exists a finite \'etale morphism 
$\mathfrak{U}_1 \to \mathfrak{U}$ 
with $\mathfrak{U}_1$ irreducible and normal
that satisfies the 
following condition:
Let
$f_1 : \mathscr{A}_1
\to \mathfrak{U}_1$ and $\mathscr{D}_1$
be the base-change of $f$ 
and $\mathscr{D}$
by this 
$\mathfrak{U}_1 \to \mathfrak{U}$,
respectively;
then 
there exists a 
finite dimensional $k$-linear subspace $\mathscr{V}$ of
$H^{0}
\left( 
\mathscr{A}_1
,
\OO_{\mathscr{A}_1} \left( 2 \mathscr{D}_1
\right)
\right)
$ 
such that $\mathscr{V}$ 
is base-point free,
that is,
the evaluation homomorphism
$\mathscr{V} \otimes_k \OO_{\mathscr{A}_1} \to \OO_{\mathscr{A}_1} 
\left( 2 \mathscr{D}_1
\right)$ is surjective.
\end{Proposition}

\Proof
Let $l$ be a prime number
with $l \neq \ch (k)$.
By Lemma~\ref{lem:bpf1},
there exist a natural number
$m$ and
a finite \'etale morphism
$\mathfrak{U}_1 \to \mathfrak{U}$
with $\mathfrak{U}_1$ irreducible
such that
\addtocounter{Claim}{1}
\begin{align} \label{eq:bpfree}
\bigcap_{\tilde{\tau} \in \mathscr{A}_1 [ l^{m} ] ( \mathfrak{U}_1 )}
\Supp
\left( T_{\tilde{\tau} }^{\ast} ( \mathscr{D}_1 )
+
T_{- \tilde{\tau}}^{\ast} ( \mathscr{D}_1)
\right)
=
\emptyset
,
\end{align}
where
$f_1 : \mathscr{A}_1 \to \mathfrak{U}_1$ and $\mathscr{D}_1$
are the base-change of $f$ and $\mathscr{D}$
by this $\mathfrak{U}_1 \to \mathfrak{U}$,
respectively.
Note that $\mathfrak{U}_1$ is normal as well as $\mathfrak{U}$.
Indeed, since the fiber of $\mathfrak{U}_1 \to \mathfrak{U}$
over any point of $\mathfrak{U}$ is a finite reduced scheme,
that follows from
\cite[Corollary of Theorem~23.9]{matsumura}.
Let $A$ be the geometric generic fiber of $f_1$
and let $D$ be the restriction of $\mathscr{D}_1$ to $A$.
Since $\ch (k)
\neq l$,
$\bigcup_{n \in \NN} A [l^{n} ]$ is dense in $A$,
and hence
there exist an $n_0 \in \NN$ 
and $\sigma \in A [l^{n_0} ]$ such that
\addtocounter{Claim}{1}
\begin{align} \label{eq:takesigma}
\sigma 
\notin
\bigcup_{\tau
\in A[ l^{m}]}
\Supp
\left(
T_{\tau }^{\ast} ( D )
\right)
.
\end{align}
Since $\ch (k)
\neq l$,
there exists a finite \'etale morphism $\mathfrak{U}_1' \to \mathfrak{U}_1$
with $\mathfrak{U}_1'$ irreducible
such that $\sigma$ extends to a section 
$\tilde{\sigma}$
of the base-change $\mathscr{A}_1'
\to 
\mathfrak{U}_1'$ of $f_1$ by this $\mathfrak{U}_1' \to \mathfrak{U}_1$.
Here, the function field of $\mathfrak{U}_1'$ is regarded
as a subfield of $\overline{K}$.
Replacing $\mathfrak{U}_1$ with $\mathfrak{U}_1'$,
we may and do assume $\mathfrak{U}_1' = \mathfrak{U}_1$.

We note by (\ref{eq:takesigma}) that
for any $\tilde{\tau}
\in \mathscr{A}_1 [ l^{m}] ( \mathfrak{U}_1 )$,
$
\tilde{\sigma} ( \mathfrak{U}_1)
\nsubseteq
\Supp
\left( T_{\tilde{\tau} }^{\ast} ( \mathscr{D}_1 )
\right)
$ holds,
and hence
the effective Cartier divisor
$\tilde{\sigma}^{\ast} \left( T_{\tilde{\tau} }^{\ast} ( \mathscr{D}_1 )
\right)$
is well-defined.

\begin{Claim} \label{cl:trivial-zerosection}
Let $\tilde{\tau} \in \mathscr{A}_1 [l^{m}] ( \mathfrak{U}_1)$.
Then
the 
Cartier divisor
$\tilde{\sigma}^{\ast} \left( T_{\tilde{\tau} }^{\ast} ( \mathscr{D}_1 )
\right)$
on $\mathfrak{U}_1$ is
trivial.
\end{Claim}

\Proof
Remark that
$\tilde{\sigma}^{\ast} \left( T_{\tilde{\tau} }^{\ast} ( \mathscr{D}_1 )
\right)
= 
(\tilde{\sigma} + \tilde{\tau})^{\ast}
( \mathscr{D}_1 )$
as Cartier divisors.
Since
$
(\tilde{\sigma} + \tilde{\tau})^{\ast} 
\left(  \mathscr{D}_1 \right) 
=
(f_1)_{\ast}
\left(
\mathscr{D}_1 \cap ( \tilde{\sigma} + 
\tilde{\tau}) ( \mathfrak{U}_1  )
\right)
$
as cycles on $\mathfrak{U}_1$,
we then find
\[
\tilde{\sigma}^{\ast} \left( T_{\tilde{\tau} }^{\ast} ( \mathscr{D}_1 )
\right)
=
(f_1)_{\ast}
\left(
\mathscr{D}_1 \cap ( \tilde{\sigma} + 
\tilde{\tau}) ( \mathfrak{U}_1  )
\right)
.
\]

Let $K_1$ be the function field of $\mathfrak{U}_1$.
Note that $K_1$ is a finite extension of $K'$.
Let $\mathfrak{B}_1 \to \mathfrak{B}'$ 
be the normalization of $\mathfrak{B}'$ in 
$K_1$.
Since $\mathfrak{U}_1$ is normal
and $\mathfrak{U}_1 \to \mathfrak{U}$
is finite,
we have
$\mathfrak{U}_1 \subset \mathfrak{B}_1$, and
since $\codim \left(
\mathfrak{B}' \setminus \mathfrak{U} , \mathfrak{B}'
\right)
\geq 2$ and $\mathfrak{U}_1 \to  \mathfrak{U}$ is finite,
we have
$\codim \left(
\mathfrak{B}_1 \setminus \mathfrak{U}_1 , \mathfrak{B}_1
\right)
\geq 2$.

Note that $A$ is nowhere degenerate by Remark~\ref{remark:new:nodeg-model}.
Let $\mathcal{H}_1$ be the pull-back of 
the line bundle $\mathcal{H}$ on $\mathfrak{B}$ 
by the composite
$\mathfrak{B}_1 \to \mathfrak{B}' \to \mathfrak{B}$.
Take a Weil divisor
$\mathfrak{D}_1$ on $\mathfrak{B}_1$
with $\mathfrak{D}_1 \cap \mathfrak{U}_1 = (f_1)_{\ast}
\left(
\mathscr{D}_1 \cap (\tilde{\sigma} + 
\tilde{\tau}) ( \mathfrak{U}_1  )
\right)$.
Then we note
\[
[ \mathfrak{D}_1 \cap \mathfrak{U}_1 ] = (f_1)_{\ast}
\left(
\cherncl_1 
(\OO_{\mathscr{A}_1} (\mathscr{D}_1)) \cdot 
[ (\tilde{\sigma} + \tilde{\tau}) ( \mathfrak{U}_1  )]
\right)
\]
as cycle classes, in particular.
Since
$L := \OO_{\mathscr{A}_1} ( \mathscr{D}_1 )|_{A}$
is even and
$0_{f_1}^{\ast} \left( \OO_{\mathscr{A}_1} ( \mathscr{D}_1 ) \right) \cong
\OO_{\mathfrak{U}_1}$
by assumption,
where $0_{f_1}$ is the zero-section of $f_1$,
it follows from
Lemma~\ref{lem:height-intersection1} that
\[
\widehat{h}_{L}
( \sigma + \tau)
=
\frac{\deg_{\mathcal{H}_1} [ \mathfrak{D}_1] }{[K_1 : K]}
.
\]
Since $\sigma + \tau$ is a torsion point,
it follows 
from the above equality
that
$\deg_{\mathcal{H}_1} [ \mathfrak{D}_1] = 0$.
Since $\mathfrak{D}_1$ is effective and $\mathcal{H}_1$
is ample, this means that $\mathfrak{D}_1$ is trivial.
Therefore 
\[
\tilde{\sigma}^{\ast} \left( T_{\tilde{\tau} }^{\ast} ( \mathscr{D}_1 )
\right)
=
(f_1)_{\ast}
\left(
\mathscr{D}_1 \cap (\tilde{\sigma} + \tilde{\tau}) ( \mathfrak{U}_1  )
\right)
=
\mathfrak{D}_1 \cap \mathfrak{U}_1
,
\]
is trivial.
Thus the claim holds.
\QED

Let us prove that for any $\tilde{\tau} \in \mathscr{A}_1
[l^{m} ] ( \mathfrak{U}_1 )$,
we have $T_{\tilde{\tau} }^{\ast} ( \mathscr{D}_1 )
+
T_{- \tau}^{\ast} ( \mathscr{D}_1)
\sim
2 \mathscr{D}_1$.
Put 
\[
\mathscr{N}
:=
\OO_{\mathscr{A}_1}
\left(
T_{\tilde{\tau} }^{\ast} ( \mathscr{D}_1 )
+
T_{- \tilde{\tau}}^{\ast} ( \mathscr{D}_1)
-
2 \mathscr{D}_1
\right)
.
\]
We then have
$\mathscr{N}|_{A} =
\OO_{A} 
\left( T_{\tau }^{\ast} ( D )
+
T_{- \tau}^{\ast} ( D )
-
2D
\right)$,
where $\tau$ is the restriction of $\tilde{\tau}$
to $A$.
By the theorem of the square
(cf. \cite[II 6 Corollary~4]{mumford}),
we have $\mathscr{N} |_A \cong \OO_A$.
This means that there exists a line bundle $\mathcal{M}$
on $\mathfrak{U}_1$ such that $\mathscr{N} = f_1^{\ast}
(\mathcal{M})$;
indeed, since $f_1$ has irreducible and reduced fibers,
if we set $\mathcal{M} := 0_{f_1}^{\ast} ( \mathscr{N})$,
then 
$\mathscr{N} = f_1^{\ast}
(\mathcal{M})$
(cf. the description of Picard functor
in the first paragraph of \cite[Chapter~0, \S~5, d)]{MFK}).
By
$0_{f_1}^{\ast} \left( \OO_{\mathscr{A}_1} ( \mathscr{D}_1 ) \right) \cong
\OO_{\mathfrak{U}_1}$
and
Claim~\ref{cl:trivial-zerosection},
we see that
$\sigma^{\ast} \left( \mathscr{N} \right) $
is trivial. Hence $\mathcal{M} 
=
\sigma^{\ast} \left( \mathscr{N} \right)
\cong \OO_{\mathfrak{U}_1}$,
and thus
$\mathscr{N} = f_1^{\ast} ( \mathcal{M} )$ is trivial.
This shows that
we have $T_{\tilde{\tau} }^{\ast} ( \mathscr{D}_1 )
+
T_{- \tilde{\tau}}^{\ast} ( \mathscr{D}_1)
\sim
2 \mathscr{D}_1$
for any
$\tilde{\tau} \in \mathscr{A}_1
[l^{m} ] ( \mathfrak{U}_1 )$.

Therefore, for each $\tilde{\tau} \in \mathscr{A}_1 [l^{m}] 
( \mathfrak{U}_1)$,
there exists a section $s_{\tilde{\tau}} \in H^{0} \left(
\mathscr{A}_1
,
\OO_{\mathscr{A}_1}
( 2 \mathscr{D}_1 )
\right)$
such that $\zero (s_{\tilde{\tau}})
= T_{\tilde{\tau} }^{\ast} ( \mathscr{D}_1 )
+
T_{- \tilde{\tau}}^{\ast} ( \mathscr{D}_1)$.
Let $\mathscr{V}$ be the $k$-vector subspace of $H^{0}
\left( 
\mathscr{A}_1
,
\OO_{\mathscr{A}_1}
( 2 \mathscr{D}_1 )
\right)$
spanned by $\{ s_{\tilde{\tau}} \}_{\tilde{\tau} \in 
\mathscr{A}_1 [l^{m}] 
( \mathfrak{U}_1)}$.
Then
by (\ref{eq:bpfree}), $\mathscr{V}$ does not have base-points.
Thus we obtain the proposition.
\QED

\subsection{Positivity of the canonical height of an ample divisor}

The purpose of this subsection is to show 
Proposition~\ref{prop:main1},
the positivity of the canonical height of an ample divisor.
To do that, we show
the following lemma,
which will be the key to the proof of the proposition.

\begin{Lemma} \label{lem:main}
Let $A$ be a nowhere degenerate abelian variety over $\overline{K}$.
Let $K'$, $\mathfrak{B}'$,
$\mathfrak{U}$,
$f : \mathscr{A} \to \mathfrak{U}$, and $0_f$
be as in Proposition~\ref{prop:basepointfree}.
Assume that $f$ is a model of $A$,
that is, $f$ has geometric generic fiber $A$.
Let $\mathscr{D}$ be an effective Cartier divisor 
on $\mathscr{A}$
that is flat 
over 
any point of $\mathfrak{U}$ of codimension $1$.
Let $D$ be the restriction of $\mathscr{D}$ to $A$.
Assume that $\OO_A(D)$ is even and ample.
Further,
suppose that
$0_f^{\ast} \left( \OO_{\mathscr{A}} (\mathscr{D}) \right) 
\cong \OO_{\mathfrak{U}}$.
Then
$A$ has non-trivial $\overline{K}/k$-trace.
\end{Lemma}

\Proof
Replacing $\mathfrak{U}$ with its open subset
whose complement in $\mathfrak{U}$ has codimension at least $2$
if necessary,
we may and do assume that $\mathscr{D}$ is flat over $\mathfrak{U}$.
Since 
$\OO_A(D)$ is even and
$0_f^{\ast} \left( \OO_{\mathscr{A}} (\mathscr{D}) \right) 
\cong \OO_{\mathfrak{U}}$,
Proposition~\ref{prop:basepointfree}
gives us
a finite surjective morphism $\mathfrak{U}_1 \to \mathfrak{U}$
with $\mathfrak{U}_1$ irreducible and normal
such that,
if $f_1 : \mathscr{A}_1 \to \mathfrak{U}_1$
and $\mathscr{D}_1$
are the base-change of $f$ and $\mathscr{D}$ 
by the morphism 
$\mathfrak{U}_1 \to \mathfrak{U}$,
respectively,
then
there exists a finite dimensional linear subspace $\mathscr{V}$
of 
$H^{0} \left( \mathscr{A}_1,
\OO_{\mathscr{A}_1} (2 \mathscr{D}_1) 
\right)$ 
such that the evaluation homomorphism
$\mathscr{V} \otimes_{k} \OO_{\mathscr{A}_1} \to \OO_{\mathscr{A}_1} 
(2 \mathscr{D}_1) $
is surjective.
Note that there exists a normal projective variety $\mathfrak{B}_1$
which contains $\mathfrak{U}_1$ with
$\codim \left( \mathfrak{B}_1 \setminus \mathfrak{U}_1
,
\mathfrak{B}_1 \right) \geq 2$;
Indeed, it suffices to let $\mathfrak{B}_1$ be
the normalization
of $\mathfrak{B}'$ in the function field of $\mathfrak{U}_1$;
see the proof of Claim~\ref{cl:trivial-zerosection} for more details.
Thus,
replacing $f : \mathscr{A} \to \mathfrak{U}$
and $\mathscr{D}$ with $f_1 : \mathscr{A}_1 \to \mathfrak{U}_1$
and $\mathscr{D}_1$, respectively,
we may and do assume that
there exists a finite dimensional linear subspace 
$\mathscr{V}$
of 
$H^{0} \left( \mathscr{A}
,
\OO_{\mathscr{A}} (2 \mathscr{D}) 
\right)$ 
such that the evaluation homomorphism
$\mathscr{V} \otimes_{k} \OO_{\mathscr{A}} \to \OO_{\mathscr{A}} (2 \mathscr{D}) $
is surjective.

To ease notation, we put $\mathscr{L} :=
\OO_{\mathscr{A}} (2 \mathscr{D})$.
Note that $0_f^{\ast} ( \mathscr{L} ) 
= 0_f^{\ast} \left( \OO_{\mathscr{A}} ( \mathscr{D} )
\right)^{\otimes 2} \cong \OO_{\mathfrak{U}}$.
Let
$
\varphi : 
\mathscr{A} \to \PP^{N}_{k}
$
be the morphism
associated to the evaluation homomorphism
$\mathscr{V} \otimes_{k} \OO_{\mathscr{A}} \to \mathscr{L} $,
where $N := \dim (\mathscr{V}) - 1$,
and let $\OO (1)$ denote the tautological line bundle on $\PP^{N}_{k}$.
Take any $\sigma \in A \left( 
\overline{K} \right)$ and let $\Delta_{\sigma}$ be the
closure of $\sigma$ in $\mathscr{A}$.
Recall that $\mathcal{H}$ is an ample line bundle on $\mathfrak{B}$,
and
let $\mathcal{H}'$ be the pull-back of $\mathcal{H}$ by
the finite morphism $\mathfrak{B}' \to \mathfrak{B}$.
Let $m$ be a positive integer such that $(\mathcal{H}')^{\otimes m}$
is very ample.
First, we apply
Proposition~\ref{prop:bertini5}~(1)
to $\mathfrak{B}'$, $\mathfrak{U}$
and $(\mathcal{H}')^{\otimes m}$.
Then
there exists a dense open subset
$V(k) \subset \left| (\mathcal{H}')^{\otimes m} \right|^{b-1}$
such that for any $\underline{D} =
(D_1 , \ldots , D_{b-1}) \in V(k)$,
the intersection
$D_1 \cap \cdots \cap D_{b-1}$
is a connected non-singular 
projective curve contained in $\mathfrak{U}$,
where recall that $b := \dim ( \mathfrak{B} ) = \dim ( \mathfrak{B}')$.

We put $L := \rest{\mathscr{L}}{A} = \OO_{A} (2D)$.
Remark that it is even and ample 
since so is $\OO_{A} (D)$ by the assumption.

\begin{Claim} \label{cl:height-intersection-curve}
Take any $\underline{D} =
(D_1 , \ldots , D_{b-1}) \in V(k)$
and set $C := D_1 \cap \cdots \cap D_{b-1}$.
Then
\[
\left(
m^{b-1} \deg \left( f|_{\Delta_{\sigma}} \right)
[K' : K]
\right)
\widehat{h}_{L} ( \sigma )
=
\deg
\left(
\cherncl_{1} 
\left(  \OO (1) \right) \cdot 
\varphi_{\ast} \left(
\left[ \Delta_{\sigma} \cap f^{-1} (C) \right]
\right)
\right)
.
\]
Furthermore,
$\widehat{h}_{L} ( \sigma ) = 0$
if and only if $\dim \left(
\varphi
\left( \Delta_{\sigma} \cap f^{-1} (C) \right)
\right) = 0$.
\end{Claim}

\Proof
By Nagata's embedding theorem
(cf. \cite[Theorem~5.7]{vojta}),
there exists a proper morphism $\bar{f} : 
\bar{\mathscr{A}} \to \mathfrak{B}'$
such that the restriction $\bar{f}^{-1} ( \mathfrak{U})
\to \mathfrak{U}$ coincides with $f :
\mathscr{A} \to \mathfrak{U}$.
Furthermore,
we can take $\bar{f}$
in such a way that
$\varphi : \mathscr{A} \to
\PP^{N}_k$ extends to a morphism $\bar{\varphi} : 
\bar{\mathscr{A}} \to \PP^{N}_k$.
Indeed,
this condition is satisfied
if we replace $\bar{\mathscr{A}}$
with the graph of the rational map 
from $\bar{\mathscr{A}}$
to $\PP^{N}_k$
given by $\varphi$.
Set $\bar{\mathscr{L}} := \bar{\varphi}^{\ast}
\left(
\OO
(1)
\right)$.
Note $\bar{\mathscr{L}} |_{\mathscr{A}} = \mathscr{L}$.
Let $\bar{\Delta}_{\sigma}$ be the closure of $\sigma$
in $\bar{\mathscr{A}}$.

From the definition of $C$,
we see that
\begin{align*}
m^{b-1} \cherncl_{1} 
\left( \mathcal{H}' \right)^{b-1} \cdot \left[ \mathfrak{B}' \right]
=
\cherncl_{1} 
\left( (\mathcal{H}')^{\otimes m} \right)^{b-1} \cdot \left[ \mathfrak{B}' 
\right]
= [C]
\end{align*}
as $1$-cycle classes on $\mathfrak{B}'$.
Since $\bar{f}$ is proper and flat over $C$
and since $C \subset \mathfrak{U}$,
it follows that
\begin{align*}
m^{b-1} \cherncl_{1} 
\left( \bar{f}^{\ast} (\mathcal{H}') \right)^{b-1} \cdot 
\left[ \bar{\Delta}_{\sigma} \right]
=
\left[ f^{-1} (C) \cap \Delta_{\sigma} \right]
\end{align*}
as $1$-cycle classes on $\bar{\mathscr{A}}$.
Therefore,
\addtocounter{Claim}{1}
\begin{align} \label{eq:left-right}
m^{b-1}
\cherncl_{1} 
\left( \bar{\varphi}^{\ast} \left( \OO (1) \right) \right) \cdot 
\cherncl_{1} ( \bar{f}^{\ast}  (\mathcal{H}') )^{\cdot (b-1)}
\cdot
\left[ \bar{\Delta}_{\sigma} \right]
=
\cherncl_{1} 
\left( \bar{\varphi}^{\ast} \left( \OO (1) \right) \right) \cdot 
\left[ \Delta_{\sigma} \cap f^{-1} (C) \right]
\end{align}
as $0$-cycle classes on $\bar{\mathscr{A}}$.
By the projection formula,
the degree of the left-hand side
of (\ref{eq:left-right}) equals
\[
m^{b-1}
\deg_{\mathcal{H}'}
\bar{f}_{\ast}
\left(
\cherncl_{1} 
\left( \bar{\varphi}^{\ast} \left( \OO (1) \right) \right) \cdot 
\left[ \bar{\Delta}_{\sigma} \right]
\right)
.
\]
Since
$\left.
\bar{\varphi}^{\ast} \left( \OO (1) \right)
\right|_{A} = L$
and $0_f^{\ast} \left( 
\bar{\varphi}^{\ast} \left( \OO (1) \right) \right)
= 0_f^{\ast} \left( \mathscr{L} \right) \cong \OO_{\mathfrak{U}}$,
that equals 
\[
\left( m^{b-1} \deg \left( f|_{\Delta_{\sigma}} \right)
[K' : K]
\right)
\widehat{h}_{L} ( \sigma )
\]
by Lemma~\ref{lem:height-intersection1}.
On the other hand,
the right-hand side
of (\ref{eq:left-right})
equals
\[
\cherncl_{1} 
\left(  \OO (1) \right) \cdot 
\varphi_{\ast} \left(
\left[ \Delta_{\sigma} \cap f^{-1} (C) \right]
\right)
.
\]
Thus we obtain
\addtocounter{Claim}{1}
\begin{align*} 
\left(
m^{b-1} \deg \left( f|_{\Delta_{\sigma}} \right)
[K' : K]
\right)
\widehat{h}_{L} ( \sigma )
=
\deg
\left(
\cherncl_{1} 
\left(  \OO (1) \right) \cdot 
\varphi_{\ast} \left(
\left[ \Delta_{\sigma} \cap f^{-1} (C) \right]
\right)
\right)
,
\end{align*}
as required.

The last assertion of the claim
follows from the fact that $\OO (1)$
is ample and $\varphi_{\ast} \left(
\left[ \Delta_{\sigma} \cap f^{-1} (C) \right]
\right)$
is effective.
\QED

\begin{Claim} \label{cl:height0-contractedtoapoint}
We have $\widehat{h}_{L} (\sigma) = 0$
if and only if 
$\varphi ( \Delta_{\sigma} )$ consists of a single point.
\end{Claim}

\Proof
First, 
we show the ``if'' part.
Suppose that $\varphi ( \Delta_{\sigma} )$ is a singleton.
Recall that $V (k)$ is the dense open subset of 
$| (\mathcal{H}')^{\otimes m}|^{b-1}$
taken just above Claim~\ref{cl:height-intersection-curve}.
Since
$V (k) \neq \emptyset$,
we take a $(D_1 , \ldots , D_{b-1}) \in V (k)$
and set $C := D_1 \cap \cdots \cap D_{b-1}$.
Then $C$ is an irreducible curve contained in $\mathfrak{U}$.
Since $\varphi ( f^{-1} ( C ) \cap \Delta_{\sigma} ) 
\subset \varphi ( \Delta_{\sigma} ) $,
we have $\dim \left( \varphi ( f^{-1} ( C ) \cap \Delta_{\sigma} ) \right) = 0$
by the assumption.
By Claim~\ref{cl:height-intersection-curve},
we conclude $\widehat{h}_{L} (\sigma) = 0$.

Let us prove the other implication.
By Proposition~\ref{prop:bertini5}~(2),
for a general
$x \in \mathfrak{U} (k)$,
the set
$
V_x(k):=
V(k) \cap 
\left( \left| (\mathcal{H}')^{\otimes m} \right|_x^{b-1} 
\right)$
is dense in $\left| (\mathcal{H}')^{\otimes m} \right|_x^{b-1}$,
where
\[
\left| (\mathcal{H}')^{\otimes m} \right|_x^{b-1} 
:=
\left\{ 
\left.
(D_1 , \ldots , D_{b-1})
\in 
\left| (\mathcal{H}')^{\otimes m} \right|^{b-1}
\ 
\right|
x \in D_1 \cap \cdots \cap D_{b-1}
\right\}
.
\]
Since $f|_{\Delta_{\sigma}} : \Delta_{\sigma}
\to \mathfrak{U}$
is generically finite and generically flat,
we take such an $x$ so that $f$ is finite and flat
over an open neighborhood of $x$.

We further set
\[
\mathscr{S} (x)
:=
\{ D_1 \cap \cdots \cap D_{b-1}
\mid
( D_1 , \ldots , D_{b-1}) \in V_x (k)
\}
.
\]
Let $y$ be a point of $\Delta_{\sigma}$ with $f (y) = x$
and
set
$\mathscr{T}_{\sigma} (y)$
to be
\begin{align*}
\{
C' \mid
\text{$C'$ is an irreducible 
component of 
$f^{-1} (C) 
\cap \Delta_{\sigma}$ 
for some
$C \in \mathscr{S} (x)$
with $y \in C'$}
\}
.
\end{align*}
Since $f|_{\Delta_{\sigma}} : \Delta_{\sigma}
\to \mathfrak{U}$ is surjective
and
is finite flat over $x$
and since any $C \in \mathscr{S} (x)$
is irreducible,
we find that
$\{ f (C') \mid C' \in \mathscr{T}_{\sigma} (y) \} = \mathscr{S} (x)$.
Therefore
\[
f \left(
\bigcup_{C' \in \mathscr{T}_{\sigma} (y)} C'
\right) 
=
\bigcup_{C \in \mathscr{S} (x)} C
.
\]
By Proposition~\ref{prop:bertini5}~(3),
this subset
is dense in $\mathfrak{B}'$.
Since 
$\Delta_\sigma$ is irreducible and
$f|_{\Delta_{\sigma}}$ is generically finite,
it follows that $\bigcup_{C' \in \mathscr{T}_{\sigma} (y)} C'$
is dense in $\Delta_{\sigma}$.

Suppose 
 $\widehat{h}_{L} ( \sigma ) = 0$.
By
Claim~\ref{cl:height-intersection-curve},
we have
$\dim \left( \varphi \left(
 \Delta_{\sigma} \cap f^{-1} (C) 
\right)\right) = 0$
for any $C \in \mathscr{S} (x)$.
This means that
for any $C' \in \mathscr{T}_{\sigma} (y)$,
we have
$\varphi \left( C' \right) = \varphi (y)$.
It follows that
$\varphi 
\left(
\bigcup_{C' \in \mathscr{T}(y)} C'
\right) = \varphi (y)$.
Since 
$\bigcup_{C' \in \mathscr{T}_{\sigma} (y)} C'$
is dense in $\Delta_{\sigma}$,
this
concludes
$\varphi 
\left(
\Delta_{\sigma}
\right) = \varphi(y)$,
which is a singleton.
Thus we complete the proof of the claim.
\QED

Note that 
$\mathscr{L}|_{f^{-1}(u)}$ is ample
for any $u \in \mathfrak{U}$.
Indeed,
since 
$\mathscr{L}$ is ample on the generic fiber
of $f$,
it is relatively ample over a dense open subset of $\mathfrak{U}$.
By the flatness of $f$, we see that $\mathscr{L}|_{f^{-1}(u)}$
is big for any $u \in \mathfrak{U}$.
Since ampleness is same as bigness
on an abelian variety,
this means that $\mathscr{L}|_{f^{-1}(u)}$ is ample
for any $u \in \mathfrak{U}$.

\begin{Claim} \label{cl:verticalcurve}
For any closed subvariety $C \subset \mathscr{A}$
of dimension $1$ such that
$\dim f (C) = 0$,
we have $\dim (\varphi (C) ) > 0$. 
\end{Claim}

\Proof
Since $\dim f (C) = 0$,
we put $u := f (C)$, which is a closed point of $\mathfrak{U}$.
Remark that $C$ is a closed curve 
of the abelian variety $f^{-1} (u)$.
Since 
$\mathscr{L}|_{f^{-1} (u)}$ is ample,
we then have 
\[
0 <
\deg ( \cherncl_1 (\mathscr{L}) \cdot [C] )
=
\deg 
\left(
\varphi^{\ast}
\left(
\cherncl_1 \left( \OO (1) \right)
\right)
\cdot [C]
\right)
=
\deg 
\left(
\cherncl_1 \left( \OO (1) \right)
\cdot \varphi_{\ast} ([C])
\right)
,
\]
which shows $\dim (\varphi (C) ) > 0$.
\QED

Put
\[
\Gamma :=
\left\{
\gamma
\mid
\text{$\gamma$ is
an irreducible closed subset
of $\mathscr{A}$
with
$\dim ( \gamma ) = b$,
and
$f (\gamma) = \mathfrak{U}$}
\right\}
.
\]
For any $\sigma \in A \left( \overline{K} \right)$,
we find $\Delta_{\sigma}
\in \Gamma$,
and
this
assignment defines a map $
A \left( \overline{K} \right) \to \Gamma$.
Note that it is surjective.
Further, put
\[
\Gamma_{0}
:=
\left\{
\gamma \in \Gamma
\left|
\ 
\dim (\varphi ( \gamma )) = 0
\right.
\right\}
\]
and
\[
A \left( \overline{K} \right)_{0}
:=
\left\{ \sigma \in
A \left( \overline{K} \right)
\left|
\ 
\widehat{h}_L (\sigma) = 0
\right.
\right\}
.
\]
Then 
by Claim~\ref{cl:height0-contractedtoapoint},
the map $
A \left( \overline{K} \right) \to \Gamma$ induces a surjective map
$\alpha : A \left( \overline{K} \right)_{0}
\to \Gamma_0$.

Since the scheme-theoretic image $\varphi ( \mathscr{A} )$ is an integral scheme,
there exists a dense open subset
$Z \subset \varphi ( \mathscr{A} )$ such that
the restriction
$\varphi^{-1} (Z)
\to Z$ is flat.
Remark that since $\mathscr{A}$ is irreducible,
$\varphi^{-1} (z)$ is pure-dimensional for any $z \in Z$.

\begin{Claim} \label{cl:dimension=b}
Let $z \in Z (k)$
and let $\Delta$ be an irreducible component
of $\varphi^{-1} (z)$.
Then
$f|_{\Delta}
: \Delta
\to \mathfrak{U}$ is a finite surjective morphism.
\end{Claim}

\Proof
First, we prove that $\dim ( \Delta ) \geq b$.
Remark that for any $\sigma \in A \left( \overline{K} \right)_{0}$,
$\alpha ( \sigma  ) = \Delta_{\sigma}$
is a closed subset of $\mathscr{A}$.
Since $A \left( \overline{K} \right)_{0}$ is dense in $A$,
$\bigcup_{\sigma
\in A \left( 
\overline{K} \right)_{0}} 
\alpha \left( \sigma \right)$
is dense
in $\mathscr{A}$.
It follows that
there exists a $\sigma \in A \left( \overline{K} \right)_{0}$
such that the generic point of $\alpha ( \sigma )$ is in 
$\varphi^{-1} ( Z )$.
Since 
$\dim \varphi (\alpha ( \sigma ) ) = 0$,
we then have $\varphi
\left( \alpha ( \sigma ) \right) \in  Z (k)$.
Furthermore,
since
$\dim (
\alpha (\sigma ) ) = b$,
it follows that
$ \dim \varphi^{-1} \left(
\varphi (\alpha ( \sigma ) )
\right) \geq b$.
Since
$\varphi$ is flat over $Z$
and $\rest{\varphi}{\varphi^{-1} (Z)} : 
\varphi^{-1} (Z) \to Z$
is surjective,
we find $\dim
\varphi^{-1} \left( z \right) \geq b$,
and thus $\dim ( \Delta  ) \geq b$.

To conclude the claim,
since $f|_{\varphi^{-1} \left( z \right)}
: \varphi^{-1} \left( z \right)
\to \mathfrak{U}$ is proper
and any irreducible component of $
\varphi^{-1} \left( z \right)$ has dimension at least $ b$,
it suffices to show that $f|_{\varphi^{-1} \left( z \right)}$ 
is quasi-finite.
We prove this by contradiction.
Suppose that 
$f|_{\varphi^{-1} \left( z \right)}$
is not quasi-finite.
Then
there exists an irreducible reduced closed curve
$C \subset \varphi^{-1} \left( z \right)$
such that $f( C) $ consists of a single point.
By Claim~\ref{cl:verticalcurve},
we have $\dim \varphi (C) > 0$.
On the other hand,
we have $\varphi (C) \subset 
\varphi \left( 
\varphi^{-1} \left( z \right)
\right)= \{ z \}$.
That is a contradiction.
\QED

Now we finally prove the lemma.
Set $\Gamma_0^Z := \{ \gamma \in \Gamma_0 \mid 
\varphi ( \gamma ) \in Z \}$.
Then for 
any $z \in Z (k)$,
Claim~\ref{cl:dimension=b} tells us that
any irreducible component of $\varphi^{-1} (z)$
is in $\Gamma_0^Z$.
This means in particular that
the map
$\Gamma_0^{Z} \to Z(k)$
which assign to each $\gamma \in \Gamma_0^{Z}$ the point 
$\varphi ( \gamma )$ is surjective.
Since $\alpha : A \left( \overline{K} \right)_0
\to \Gamma_0$ is surjective,
we then find that
\addtocounter{Claim}{1}
\begin{align} \label{eq:cardinality1}
\# A \left( \overline{K} \right)_0
\geq \# \Gamma_0
\geq \# \Gamma_0^{Z}
\geq \# Z (k)
.
\end{align}

We set $n := \dim (A) > 0$.
Since Claim~\ref{cl:verticalcurve}
shows that the restriction of $\varphi$
to a closed fiber of $f$ is a finite morphism,
we note that $\dim (Z) = \dim \varphi ( \mathscr{A} ) \geq n$.

Now, assume that $k$ has uncountably infinite cardinality,
and
we prove that $A$ has non-trivial $\overline{K}/k$-trace
under this condition first.
Let $\aleph_0$ denote the countably infinite cardinality.
Then
since $\dim (Z) 
> 0$,
we have
$\# Z (k) > \aleph_{0}$.
It follows from inequality
(\ref{eq:cardinality1})
that
$\# A \left( \overline{K} \right)_0 > \aleph_{0}$.
Since $\# A \left( \overline{K} \right)_{tor} = \aleph_0$,
this means that $A \left( \overline{K} \right)$ 
has height $0$ points other than torsion points.
By \cite[Chapter~6, Theorem~5.5]{lang2},
we conclude that $A$ has non-trivial $\overline{K}/k$-trace
in this case.

We consider the general case.
Take an algebraically closed field extension $k'$ of $k$
such that $\# k' > \aleph_0$.
Then 
$\mathfrak{B} 
\otimes_k k'$
is a normal projective variety over $k'$.
Further,
$\mathfrak{U} \otimes_k k'$,
$f \otimes_k k' 
:
\mathscr{A} \otimes_k k'
\to \mathfrak{U} \otimes_k k'$
and $\mathscr{D} \otimes_k k'$
satisfy all the conditions of the lemma
as a variety, a morphism and a Cartier divisor over 
the constant field $k'$,
respectively.
Let $F$ be the function field of $\mathfrak{B} 
\otimes_k k'$.
Since we know that the lemma holds if the constant field
has cardinality greater than $\aleph_0$,
it turns out that $A \otimes_{\overline{K}} \overline{F}$
has non-trivial $\overline{F} / k'$-trace.
By \cite[Lemma~A.1]{yamaki7},
it follows that $A$ has non-trivial $\overline{K} / k
$-trace.
Thus we obtain the lemma.
\QED

Now we show the positivity assertion.

\begin{Proposition} \label{prop:main1}
Let $A$ be a nowhere degenerate abelian variety over $\overline{K}$
with trivial $\overline{K}/k$-trace
and let
$L$ be an even ample line bundle on $A$.
Further, let $X$ be an
effective
ample divisor on $A$.
Then we have
$\widehat{h}_{L} (X) >0$.
\end{Proposition}

\Proof
There exists a torsion point $\tau \in A \left( 
\overline{K} \right)$ such that
$0 \notin \Supp  T_{\tau}^{\ast} (X) $.
Since
being $\widehat{h}_{L} (X) >0$
is
equivalent to 
being $\widehat{h}_L (T_{\tau}^{\ast} (X) ) > 0$,
we may assume that $0 \notin X$.
We put $D := X + [-1]^{\ast} (X)$,
where ``$+$'' means the addition of the divisors.
Then $D$ is an even ample divisor with
$0 \notin \Supp (D)$.
Since whether or not $\widehat{h}_{L} (X) >0$
does not depend on the choice of even ample line bundles $L$,
we may assume that $L = \OO_A \left( D \right)$.
Furthermore,
since $[-1]$ is an automorphism and $[-1]^{\ast} (L) = L$,
we have $\widehat{h}_{L} (X) = \widehat{h}_{L} (D) / 2$,
so that it suffices to show that $\widehat{h}_{L} (D) > 0$.

By Proposition~\ref{prop:model1},
there exist
a finite extension $K'$ of $K$,
an open subset $\mathfrak{U}$ 
of 
$\mathfrak{B}'$
with
$\codim ( \mathfrak{B}' 
\setminus \mathfrak{U} , \mathfrak{B}')
\geq 2$,
where $\mathfrak{B}'$
is the normalization
of $\mathfrak{B}$
in $K'$,
an abelian scheme $f : \mathscr{A} \to \mathfrak{U}$
with zero-section $0_f$,
and 
a line bundle $\mathscr{L}$ on $\mathscr{A}$
such that
$\left( f , \mathscr{L} \right)$ is a model of
$(A ,L)$
over $\mathfrak{U}$
with $0_{f}^{\ast} ( \mathscr{L} ) \cong \OO_{\mathfrak{U}}$.
Replacing this $K'$ with its finite extension if necessary,
we may and do assume that
$D$ can be defined over the function field 
$K'$.
Let 
$\mathscr{D}$ be the closure of $D$ in 
$\mathscr{A}$.
Since $D$ can be defined over $K'$,
the restriction of $\mathscr{D}$ to $A$ coincides with
$D$.
Note that $\mathscr{D}$ is flat over any 
codimension $1$ point of $\mathfrak{U}$.
Since $A$ has trivial $\overline{K}/k$-trace,
it follows from Lemma~\ref{lem:main}
that $0_{f}^{\ast} \left( \OO_{\mathscr{A}} ( \mathscr{D} ) \right)$ is
a non-trivial line bundle on $\mathfrak{U}$.
Since $0_f ( \mathfrak{U} ) \nsubseteq \Supp ( \mathscr{D} )$
by assumption,
$0_{f}^{\ast} ( \mathscr{D} )$ is a non-trivial effective Cartier divisor
on $\mathfrak{U}$.

By
Nagata's embedding theorem
(cf. \cite[Theorem~5.9]{vojta}),
there exists a proper morphism $\bar{f} : 
\bar{\mathscr{A}} \to \mathfrak{B}'$
such that
$\bar{f}$ equals $f$ over $\mathfrak{U}$
and that
the Cartier divisors $\mathscr{D}$
and $f^{\ast} \left( 0_f^{\ast} ( \mathscr{D} ) \right)$ 
on $\mathscr{A}$
extends to  Cartier divisors $\bar{\mathscr{D}}$
and $\bar{\mathscr{E}}$
on 
$\bar{\mathscr{A}}$
respectively.
Since
$\mathscr{D}$ is flat over any point of $\mathfrak{U}$
of codimension $1$ and since
 $\codim \left( \mathfrak{B}' \setminus  
\mathfrak{U} , \mathfrak{B}'
\right) \geq 2$,
$\bar{\mathscr{D}}$ 
is flat over any point of $\mathfrak{B}'$ of codimension $1$.
We set 
$\bar{\mathscr{L}} := 
\OO_{\bar{\mathscr{A}}} 
\left( \bar{\mathscr{D}} 
\right)$,
$\bar{\mathscr{N}} := 
\OO_{\bar{\mathscr{A}}} 
\left( \bar{\mathscr{E}} \right)$
and
$\bar{\mathscr{L}}_0 :=
\bar{\mathscr{L}} \otimes 
\bar{\mathscr{N}}^{ \otimes (-1)}$.

Let $\mathcal{H}'$ be the pull-back of $\mathcal{H}$
by the finite morphism $\mathfrak{B}'
\to \mathfrak{B}$
and
let $m'$ be an integer such that $\mathcal{H}''
:=
( \mathcal{H}' )^{\otimes m'}$ is very ample.
Put $n := \dim (A)$.
Then, $\widehat{h}_L (D) > 0$
if and only if
\addtocounter{Claim}{1}
\begin{align} \label{eq:goalforpositive}
\deg
\left(
\cherncl_1 \left( \bar{\mathscr{L}}_0 \right)^{\cdot n}
\cdot
\cherncl_1 \left( \bar{\mathscr{N}} \right)
\cdot
\bar{f}^{\ast}
\cherncl_1
( \mathcal{H}'')^{\cdot (b-1)}
\cdot
\left[ \bar{\mathscr{A}} \right]
\right)
>0
.
\end{align}
Indeed,
since
$\mathscr{L}_0 |_{A} = L$,
$0_f^{\ast} \left( \bar{\mathscr{L}}_0 \right) \cong
\OO_{\mathfrak{U}}$
and since
$\bar{\mathscr{D}}$ is flat over any point of $\mathfrak{B}'$
of codimension $1$,
we have,
by \cite[Lemma~2.6]{yamaki7},
\begin{align*}
\widehat{h}_{L} ( D )
=
\frac{
\deg
\left(
\cherncl_1
(\mathcal{H}'')^{\cdot (b-1)}
\cdot
\bar{f}_{\ast}
\left(
\cherncl_1 ( \bar{\mathscr{L}}_0 )^{\cdot n}
\cdot
\left[  \bar{\mathscr{D}} \right]
\right)
\right)
}{(m')^{b-1} [K':K]}
.
\end{align*}
By Remark~\ref{rem:heightabelianvariety},
we have
\begin{align*} 
\deg
\left(
\cherncl_1
( \mathcal{H}'')^{\cdot (b-1)}
\cdot
\bar{f}_{\ast}
\left(
\cherncl_1 \left( \bar{\mathscr{L}}_0 \right)^{\cdot (n+1)}
\cdot
\left[ \bar{\mathscr{A}} \right]
\right)
\right)
=
0
.
\end{align*}
Since
$\OO_{\bar{\mathscr{A}}} \left( \bar{\mathscr{D}} \right) =
\bar{\mathscr{L}}_0 \otimes \bar{\mathscr{N}}$,
it follows that
\begin{align*}
&\deg
\left(
\cherncl_1
( \mathcal{H}'')^{\cdot (b-1)}
\cdot
\bar{f}_{\ast}
\left(
\cherncl_1 \left( \bar{\mathscr{L}}_0 \right)^{\cdot n}
\cdot
\left[  \bar{\mathscr{D}} \right]
\right)
\right)
\\
&=
\deg
\left(
\cherncl_1
( \mathcal{H}'')^{\cdot (b-1)}
\cdot
\bar{f}_{\ast}
\left(
\cherncl_1 \left( \bar{\mathscr{L}}_0 \right)^{\cdot n}
\cdot
\left(
\cherncl_1 \left( \bar{\mathscr{L}}_0 \right) +
\cherncl_1 \left( \bar{\mathscr{N}} \right)
\right)
\cdot
\left[ \bar{\mathscr{A}} \right]
\right)
\right)
\\
&=
\deg
\left(
\cherncl_1
( \mathcal{H}'')^{\cdot (b-1)}
\cdot
\bar{f}_{\ast}
\left(
\cherncl_1 \left( \bar{\mathscr{L}}_0 \right)^{\cdot n}
\cdot
\cherncl_1 \left( \bar{\mathscr{N}} \right)
\cdot
\left[ \bar{\mathscr{A}} \right]
\right)
\right)
\\
&=
\deg
\left(
\cherncl_1 \left( \bar{\mathscr{L}}_0 \right)^{\cdot n}
\cdot
\cherncl_1 \left( \bar{\mathscr{N}} \right)
\cdot
\bar{f}^{\ast}
\cherncl_1
( \mathcal{H}'')^{\cdot (b-1)}
\cdot
\left[ \bar{\mathscr{A}} \right]
\right)
.
\end{align*}
Thus $\widehat{h}_{L} ( D) >0$
if and only if (\ref{eq:goalforpositive}) holds.

Let us prove (\ref{eq:goalforpositive}).
Since $\codim \left(
\mathfrak{B}' \setminus \mathfrak{U},
\mathfrak{B}'
\right) \geq 2$
and $\mathcal{H}''$ is a very ample line bundle on $\mathfrak{B}'$,
there exists a proper curve $C$ on $\mathfrak{B}'$
such that $C \subset \mathfrak{U}'$,
$C$ intersects with $0_f^{\ast} ( \mathscr{D} )$ properly,
and $\cherncl_1 ( \mathcal{H}'')^{\cdot (b-1)} \cdot [ \mathfrak{B}' ]
= [ C ]$
as cycle classes on $\mathfrak{B}'$
(cf. Proposition~\ref{prop:bertini5}).
Set $Z := C \cap 0_f^{\ast} ( \mathscr{D} )$.
Since $\mathcal{H}''$ is ample
and $0_f^{\ast} ( \mathscr{D} )$
is a non-trivial effective Cartier divisor,
we note
that $Z$ is a non-trivial
effective $0$-cycle, and hence 
$\deg [ Z ] = \mathrm{length} (\OO_Z) > 0$.
Since $\bar{f}$ is flat over $\mathfrak{U}$,
we note
$\bar{f}^{\ast}
\cherncl_1
( \mathcal{H}'')^{\cdot (b-1)}
\cdot
\left[ \bar{\mathscr{A}} \right]
= \left[ \bar{f}^{-1} ( C ) \right]$.
Since $\bar{f}^{-1} (C) \subset \mathscr{A}$
and
$\bar{\mathscr{N}} |_{\mathscr{A}} = f^{\ast}
\left(
\OO_{\mathfrak{U}}
\left(
0_{f}^{\ast} ( \mathscr{D} )
\right)
\right)$,
we find that
$\cherncl_{1} \left(
\bar{\mathscr{N}}
\right) \cdot \left[ \bar{f}^{-1} ( C ) \right] 
= \bar{f}^{\ast}  \left[
0_{f}^{\ast} ( \mathscr{D} ) \cap C
\right]
=
\bar{f}^{\ast} [ Z ]$.
It follows that
\addtocounter{Claim}{1}
\begin{align} \label{eq:degof0cycle:new13}
\deg
\left(
\cherncl_1 \left( \bar{\mathscr{L}}_0 \right)^{\cdot n}
\cdot
\cherncl_1 \left( \bar{\mathscr{N}} \right)
\cdot
\bar{f}^{\ast}
\cherncl_1
( \mathcal{H}'')^{\cdot (b-1)}
\cdot
\left[ \bar{\mathscr{A}} \right]
\right)
=
\deg
\left(
\cherncl_1 \left( \bar{\mathscr{L}}_0 \right)^{\cdot n}
\cdot
\bar{f}^{\ast} [ Z ]
\right)
.
\end{align}

Remark that
since 
$\bar{f}$ is flat over $\mathfrak{U}$,
for any $u , u'\in \mathfrak{U}$,
we have
\addtocounter{Claim}{1}
\begin{align} \label{align:comparison-closed-generic}
\deg
\left(
\cherncl_1 \left( \rest{\bar{\mathscr{L}}_0}{\bar{f}^{-1} (u)} \right)^{\cdot n}
\cdot
\left[ \bar{f}^{-1} (u) \right]
\right)
=
\deg
\left(
\cherncl_1 \left( \rest{\bar{\mathscr{L}}_0}{\bar{f}^{-1} (u')} \right)^{\cdot n}
\cdot
\left[ \bar{f}^{-1} (u') \right]
\right)
.
\end{align}
Here, suppose that 
$u \in \mathfrak{U} (k)$. Then the left-hand side 
of (\ref{align:comparison-closed-generic}) 
equals
$\deg
\left(
\cherncl_1 \left( \bar{\mathscr{L}}_0 \right)^{\cdot n}
\cdot
\left[ \bar{f}^{-1} (u) \right]
\right)$.
On the other hand, suppose that $u'$ is the generic point
of $\mathfrak{U}$.
Then
the right-hand side of (\ref{align:comparison-closed-generic})
equals
$\deg
\left(
\cherncl_1 \left( \left. \bar{\mathscr{L}}_0 \right|_{A}
 \right)^{\cdot n} \cdot [ A ]
\right)$.
Thus we have
\addtocounter{Claim}{1}
\begin{align} \label{align:forreferee}
\deg
\left(
\cherncl_1 \left( \bar{\mathscr{L}}_0 \right)^{\cdot n}
\cdot
\left[ \bar{f}^{-1} (u) \right]
\right)
=
\deg
\left(
\cherncl_1 \left( \rest{\bar{\mathscr{L}}_0}{A} \right)^{\cdot n}
\cdot
\left[ A \right]
\right)
.
\end{align}
Since $Z$ is a $0$-cycle,
this is a sum of closed points.
Further, the support of $Z$
is contained in $\mathfrak{U}$.
Therefore, it follows form
(\ref{align:forreferee})
that
\addtocounter{Claim}{1}
\begin{align} \label{align:new1313}
\deg
\left(
\cherncl_1 \left( \bar{\mathscr{L}}_0 \right)^{\cdot n}
\cdot
f^{\ast} [Z]
\right)
=
 \deg [Z] 
\cdot
\deg
\left(
\cherncl_1 \left( \left. \bar{\mathscr{L}}_0 \right|_{A}
 \right)^{\cdot n} \cdot [ A ]
\right)
.
\end{align}
Since $\left. \bar{\mathscr{L}}_0 \right|_{A} = L$ is ample,
$\deg
\left(
\cherncl_1 \left( \left. \bar{\mathscr{L}}_0 \right|_{A}
 \right)^{\cdot n} \cdot [ A ]
\right) > 0$.
Thus by (\ref{eq:degof0cycle:new13}) and
(\ref{align:new1313}), we obtain
\[
\deg
\left(
\cherncl_1 \left( \bar{\mathscr{L}}_0 \right)^{\cdot n}
\cdot
\cherncl_1 \left( \bar{\mathscr{N}} \right)
\cdot
\bar{f}^{\ast}
\cherncl_1
( \mathcal{H}'')^{\cdot (b-1)}
\cdot
\left[ \bar{\mathscr{A}} \right]
\right)
>0
.
\]
This completes
the proof of the theorem.
\QED

This is a small remark:
The assumption that $X$ is ample
is made
in Proposition~\ref{prop:main1},
but in fact
it can be removed;
see 
Proposition~\ref{prop:nondegtrivialtrace}
and
Remark~\ref{rem:main1-nondegtrivialtrace}.

\section{Non-density of small points} \label{sect:GBC}

\subsection{Non-density of small points on closed subvarieties}

In this section, we focus on 
the non-density of small points on
specified closed subvarieties
of a given abelian variety:
Suppose we are given an abelian 
variety $A$ over $\overline{K}$;
let $X$ be a closed subvariety;
then we consider the following assertion for $X$.

\begin{Conjecture} \label{conj:WGBC}
If
$X$ has dense small points,
then $X$ is a special subvariety.
\end{Conjecture}

We regard Conjecture~\ref{conj:WGBC}
as a conjecture
for $X$
while
we regard
the geometric Bogomolov conjecture 
(Conjecture~\ref{conj:GBCforAVint})
as a conjecture
for $A$.
The geometric Bogomolov conjecture for $A$
is equivalent to Conjecture~\ref{conj:WGBC}
for all closed subvarieties $X$ of $A$.

\begin{Remark} \label{rem:trivialcase}
For a closed subvariety of dimension
$\dim (A)$,
Conjecture~\ref{conj:WGBC} holds trivially.
For a closed subvariety of dimension $0$,
it is classically known that Conjecture~\ref{conj:WGBC} holds
(cf. \cite[Remark~7.4]{yamaki6}).
Therefore,
the geometric Bogomolov conjecture holds for abelian varieties
of dimension at most $1$.
\end{Remark}

\subsection{Relative height} \label{subsect:relativeheight}

In this subsection,
we fix the following.
Let $A$ be a nowhere degenerate abelian variety over $\overline{K}$
and
let
$L$ be an even ample line bundle on $A$.
Let $\widetilde{Y}$ be a variety over $k$
and set $Y := \widetilde{Y} \otimes_k \overline{K}$.
Let $p : Y \times A \to Y$ be the natural projection,
which is an abelian scheme over $Y$.
Let $X$ be a closed subvariety of $Y \times A$
such that $p (X) = Y$.

In \cite[\S~4.2]{yamaki7},
we defined an
$\RR$-valued function
$\mathbf{h}_{X/Y}^L$
over some 
subset 
$\widetilde{Y}_{\mathrm{pd}}$
of $\widetilde{Y}$,
called the relative height function.
Roughly speaking,
it is a function 
that assigns to each $\tilde{y} \in \widetilde{Y}_{\mathrm{pd}}$
the canonical height of the fiber of $\rest{p}{X}$ over
a geometric point of $Y$
which corresponds to $\tilde{y}$ in a canonical way.
We defined indeed $\mathbf{h}_{X/Y}^L$
not only for closed points of $\tilde{Y}$
but also for general points.
However, in the later argument in
this paper, it will be enough for us to
consider 
$\mathbf{h}_{X/Y}^L$ only over a set of closed points.

Therefore, we recall how
$\mathbf{h}_{X/Y}^L$ is given for closed points.
Let $\tilde{y}$ be a point in $\widetilde{Y} (k)$.
Note that $\tilde{y}$ is regarded as a point of $Y \left(
\overline{K}
\right)$ naturally.
Indeed, since $\tilde{y}$ 
is a morphism
$\Spec (k) \to \widetilde{Y}$,
taking the fiber product 
with $\Spec \left( \overline{K} \right)$
over $\Spec (k)$,
we obtain $\Spec \left( \overline{K} \right) \to \widetilde{Y} \otimes_k
\overline{K} = Y$,
which is
a point in $Y \left(
\overline{K}
\right)$.
Thus we have $\widetilde{Y}(k) \hookrightarrow Y 
\left( \overline{K}
\right)$.
We denote
by $\tilde{y}_{\overline{K}}$
the point in $Y \left(
\overline{K}
\right)$ 
corresponding
to $\tilde{y}$
via this inclusion.

The fiber $p^{-1} \left( \tilde{y}_{\overline{K}} \right)$ is an abelian
variety; indeed,
the second projection $Y \times A \to A$ restricts
to an isomorphism
$p^{-1} \left( \tilde{y}_{\overline{K}} \right) \cong A$.
We set
$X_{\tilde{y}_{\overline{K}}}
:= \rest{p}{X}^{-1} \left( \tilde{y}_{\overline{K}} \right)
= p^{-1} \left(
\tilde{y}_{\overline{K}}
\right) \cap X$,
which
is a closed subscheme of that abelian variety.
If follows that
if $X_{\tilde{y}_{\overline{K}}}$ has pure dimension,
then the canonical height 
$\widehat{h}_{L} \left( X_{\tilde{y}_{\overline{K}}} \right)$
of $X_{\tilde{y}_{\overline{K}}}$
with respect to $L$ is defined.
Thus we set
\addtocounter{Claim}{1}
\begin{align} \label{align:pd}
\widetilde{Y}_{\mathrm{pd}} (k)
:=
\left\{
\widetilde{y} \in \widetilde{Y} (k)
\left|
\ 
\text{$X_{\tilde{y}_{\overline{K}}}$
has pure dimension
$\dim (X) - \dim (Y)$}
\right.
\right\}
,
\end{align}
and
by definition,
$\mathbf{h}_{X/Y}^L $
is given by
\addtocounter{Claim}{1}
\begin{align} \label{eq:defofmathbfh}
\mathbf{h}_{X/Y}^L ( \tilde{y} )
= \widehat{h}_{L} \left( X_{\tilde{y}_{\overline{K}}} \right)
\end{align}
for $\tilde{y} \in \widetilde{Y}_{\mathrm{pd}} (k)$.

\begin{Remark} \label{remark:closedpointcase}
In \cite[\S~4.2]{yamaki7},
for any $\tilde{y} \in \widetilde{Y}$, we defined
a geometric point
$\overline{\tilde{y}_K}$ of $Y$,
which coincides with above $\tilde{y}_{\overline{K}}$
when $\tilde{y} \in \widetilde{Y} (k)$.
Further, 
we defined
$Y_{\mathrm{pd}}$ to be the set of points $\tilde{y} \in 
\widetilde{Y}$
such that
$X_{\overline{\tilde{y}_K}}$ has pure dimension
$\dim (X) - \dim (Y)$,
where $X_{\overline{\tilde{y}_K}}$ is the fiber of $\rest{p}{X}
: X \to Y$ over $\overline{\tilde{y}_K}$.
Thus the above $\widetilde{Y}_{\mathrm{pd}} (k)$
actually 
equals $Y_{\mathrm{pd}} \cap \widetilde{Y} (k)$
and the relative height 
$\mathbf{h}_{X/Y}^L$
described 
in (\ref{eq:defofmathbfh})
coincides with the restriction of the relative height
in \cite[\S~4.2]{yamaki7}
to the subset
$\widetilde{Y}_{\mathrm{pd}} (k)$.
\end{Remark}

\begin{Lemma} \label{lem:relativeheight0-densesmallpoints}
Let $A$,
$L$, $\widetilde{Y}$, $p : Y \times A \to Y$,
and $X$ be as above.
Let $\widetilde{Y}_{\mathrm{pd}} (k)$ be as in 
(\ref{align:pd}).
Let $\tilde{y}$ be a point in $ \widetilde{Y}_{\mathrm{pd}} (k)$
and
let $\tilde{y}_{\overline{K}} 
\in Y \left(
\overline{K}
\right)$ be the corresponding point via the
natural inclusion $\widetilde{Y}(k) \hookrightarrow 
Y \left(
\overline{K}
\right)$.
Then the following are equivalent to each other:
\begin{enumerate}
\renewcommand{\labelenumi}{(\alph{enumi})}
\item
Any irreducible component of $X_{\tilde{y}_{\overline{K}}}
=
\left( p|_{X} \right)^{-1} (\tilde{y}_{\overline{K}})$ has 
dense small points as a closed subvariety of 
$p^{-1} 
\left(
\tilde{y}_{\overline{K}}
\right)
\cong A$;
\item
$\mathbf{h}^{L}_{X/Y} ( \tilde{y} ) = 0$.
\end{enumerate}
\end{Lemma}

\Proof
Noting Remark~\ref{remark:closedpointcase},
we 
immediately obtain the lemma from
(\ref{eq:defofmathbfh})
and \cite[Remark~4.5]{yamaki7}.
\QED

\subsection{Proof of the results}

In this subsection, we establish
Theorem~\ref{thm:divisorialcase}:
Conjecture~\ref{conj:WGBC}
holds for any closed subvariety of $A$ of codimension $1$.
Furthermore, using this theorem, 
we show that 
Conjecture~\ref{conj:WGBC}
holds for any closed subvariety of $A$ of dimension $1$
(cf. Theorem~\ref{thm:curvecase}).

Before giving the proofs of the theorems,
we prove a couple of technical lemmas.

\begin{Lemma} 
\label{lem:rev-final}
Let $A$ be a nowhere degenerate abelian variety over $\overline{K}$.
Let
$\widetilde{B}$ be an abelian variety over $k$
and set $B := \widetilde{B} \otimes_k \overline{K}$.
Let $X$ be a closed subvariety of $B \times A$.
Let $\pr_B : B \times A \to B$ be the first projection
and
set $Y := \pr_B (X)$.
Suppose that $X$ has dense small points.
Then
the following holds.
\begin{enumerate}
\item
There exists a closed subvariety $\widetilde{Y} \subset \widetilde{B}$
such that $Y = \widetilde{Y} \otimes_k \overline{K}$.
\item
There exists a dense open subset $\widetilde{V} \subset \widetilde{Y}$
with the following property:
For any $\tilde{y} \in \widetilde{V} (k)$,
each irreducible component of
$\pr_A \left( \pr_B^{-1} 
\left(
\tilde{y}_{\overline{K}}
\right) \cap X \right) \subset A$
has dense small points,
where
$\tilde{y}_{\overline{K}} \in Y \left( \overline{K} \right)$
denotes the point in $Y \left( \overline{K} \right)$ 
corresponding to $y$ via
$\widetilde{V} (k) \hookrightarrow
Y \left( \overline{K} \right)$
and
where $\pr_A : B \times A \to A$ is the second projection.
\end{enumerate}
\end{Lemma}

\Proof
Let $L$ be an even ample line bundle on $A$.
We have (1) 
by
\cite[Proposition~5.1]{yamaki7}.
Furthermore, this proposition says 
that there exists a dense open subset $\widetilde{V}
\subset
\widetilde{Y}
$
such that
$\widetilde{V} (k)
\subset
\widetilde{Y}_{\mathrm{pd}} (k)$ 
and such that
$\mathbf{h}_{X/Y}^{L} (\tilde{y}) = 0$
for any $\tilde{y} \in \widetilde{V} (k)$
(note also Remark~\ref{remark:closedpointcase}).
By Lemma~\ref{lem:relativeheight0-densesmallpoints},
it follows that for
any $\tilde{y} 
\in \widetilde{V}(k)
$,
each irreducible component of 
$\left( \pr_B|_X \right)^{-1} 
\left(
\tilde{y}_{\overline{K}}
\right) \subset \left\{ \tilde{y}_{\overline{K}} \right\} \times A$ 
has dense small points.
Since $\pr_A$
induces an isomorphism 
$\left\{ \tilde{y}_{\overline{K}} \right\} \times A \cong A$,
any irreducible component of
$\pr_A \left( \pr_B^{-1} \left(
\tilde{y}_{\overline{K}}
\right) \cap X \right) \subset A$
has dense small points.
\QED

\begin{Lemma} \label{lem:subandquotient}
Let $A$ be an abelian variety over $\overline{K}$.
Let $A_1$ be an abelian subvariety of $A$
and let $A_2$ be
a quotient abelian variety of $A$.
Then the following hold.
\begin{enumerate}
\item
If $A$ has trivial $\overline{K}/k$-trace,
then $A_1$ and $A_2$ have trivial $\overline{K}/k$-trace.
\item
If $A$ is nowhere degenerate,
then $A_1$ and $A_2$ are nowhere degenerate.
\end{enumerate}
\end{Lemma}

\Proof
The assertion (2) follows from \cite[Lemma~7.8]{yamaki6}
immediately.
To show (1), suppose that $A$ has trivial $\overline{K}/k$-trace.
Let $\widetilde{B}$ be any abelian variety over $k$
and
set $B:= \widetilde{B} \otimes_k \overline{K}$.

Let $\phi_1 : B \to A_1$ be a homomorphism.
Since $A_1 \subset A$, $\phi_1$ is then regarded as a
homomorphism from $B$ to $A$.
Since $A$ has trivial $\overline{K}/k$-trace,
this is trivial.
Thus $\phi_1$ is trivial, which shows that $A_1$ has trivial 
$\overline{K}/k$-trace.

To show that $A_2$ has trivial $\overline{K}/k$-trace,
let $\phi_2 : B \to A_2$ be a homomorphism.
By the Poincar\'e complete reducibility theorem,
there exists a finite homomorphism $\psi : A_2 \to A$.
Since $A$ has trivial $\overline{K}/k$-trace,
the composite $\psi \circ \phi_2 : B \to A$ is trivial.
Since $\psi$ is finite, it follows that $\phi_2$ is
trivial.
This completes the proof.
\QED

Let $A$ be an abelian variety $\overline{K}$.
By \cite[Lemma~7.9]{yamaki6},
there exists a unique abelian subvariety of $A$
that is
maximal among the nowhere degenerate abelian subvarieties. 
This abelian subvariety is called the
\emph{maximal nowhere degenerate abelian subvariety of $A$}
(\cite[Definition~7.10]{yamaki6}).

Now, we show that Conjecture~\ref{conj:WGBC}
holds for
codimension $1$ subvarieties.

\begin{Theorem} \label{thm:divisorialcase}
Let $A$ be an abelian variety over $\overline{K}$
and let
$X$ be a closed subvariety of $A$ of codimension $1$.
Then
if
$X$ has dense small points,
then $X$ is a special subvariety.
\end{Theorem}

\Proof
We prove the theorem in three steps.

\smallskip

\emph{Step~1. In this step, assume that  $A$ is nowhere degenerate and
has trivial $\overline{K}/k$-trace.
Then
we prove that $X$
does not have dense small points.}
(Note that then we have the theorem for such an $A$,
because there does not exist $X$ 
that satisfies the assumption of the theorem.)
Since $X$ is a non-trivial effective divisor on $A$,
the linear system $|2X|$ gives a morphism
$\varphi: A \to \PP^N$.
There exists a hyperplane $H$ in $\PP^N$
such that $\varphi^{\ast} (H) = 2X$.
Let $S (X)$ be a closed subgroup 
of $A$ given by
$
S(X) \left( \overline{K} \right) =
\left\{
a \in A \left( \overline{K} \right)
\mid
T_{a}^{\ast}
(X) = X
\right\}
$,
where $T_a$ is the translate by $a$ and
``$T_{a}^{\ast}
(X) = X$'' means an equality as divisors.
Let $S(X)^0$ be the connected component of $S(X)$ with $0
\in S(X)^0$
and let $\phi : A \to A_0 := A/S(X)^0$ be the quotient homomorphism.
By Lemma~\ref{lem:subandquotient}, 
the abelian variety $B$ is nowhere degenerate and has trivial
$\overline{K}/k$-trace.

By \cite[p.88, Remarks on effective divisors by Nori]{mumford},
$\varphi$ factors though $\phi$,
i.e.,
there exists a finite morphism $\psi :  A_0 \to \PP^N$
such that $\varphi = \psi \circ \phi$.
Set $E := \psi^{\ast} ( H )$.
Then $E$ is an effective ample divisor on $A_0$.
Furthermore, 
since $\phi$ is surjective and $2X = \varphi^{\ast} ( H )$,
we have $\phi (X) = E$ as closed subsets of $A_0$.
Since $A_0$ is nowhere degenerate and has trivial
$\overline{K}/k$-trace,
Proposition~\ref{prop:main1}
tells us that
$E$ does not have dense small points.
Therefore
by \cite[Lemma~2.1]{yamaki5}, $X$ does not have dense small points either.

\smallskip

\emph{Step~2. We prove the theorem for the case where 
$A$ is nowhere degenerate.}
Suppose here on that $X$ has dense small points.
Let 
$\left( \widetilde{A}^{\overline{K}/k} , \Tr_A \right)$
be the $\overline{K}/k$-trace of $A$.
To ease notation, we put $\widetilde{B} := \widetilde{A}^{\overline{K}/k}$
and $B := \widetilde{B}
\otimes_k \overline{K}$.
Let
$\mathfrak{t}$ be the image of 
$\Tr_A$
and set
$A_1:= A/\mathfrak{t}$.
Since $A$ is nowhere degenerate, $A_1$
is nowhere degenerate and has trivial $\overline{K}/k$-trace
(cf. \cite[Lemma~7.8~(2)]{yamaki6}
and \cite[Remark~5.4]{yamaki7}).
Then 
by the Poincar\'e
complete reducibility theorem,
$A$ is isogenous to $\mathfrak{t} \times A_1$.
Since $B$ is isogenous to $\mathfrak{t}$
(cf. \cite[Lemma~1.4]{yamaki5}),
it follows that
there exists an isogeny $\phi : A \to B \times A_1$.
Set $X' := \phi (X)$, which is a closed subvariety of $B \times A_1$
of codimension $1$.
Since $X$ has dense small points,
so does $X'$
(cf. \cite[Lemma~2.1]{yamaki5}).
By Lemma~\ref{lem:rev-final}~(1),
there exists a closed subvariety $\widetilde{Y}
\subset \widetilde{B}$ such that
$\pr_{B} (X') = Y:= \widetilde{Y} \otimes_k \overline{K}$,
where $\pr_{B} : B \times A_1 \to B$ is the canonical projection.
Let $\widetilde{Y}_{\mathrm{pd}} (k)$ be as on (\ref{align:pd}).
By Remark~\ref{remark:closedpointcase}
and
Lemma~\ref{lem:rev-final}~(2),
there exists a subset $\widetilde{S} \subset \widetilde{Y}_{\mathrm{pd}} 
(k)$
such
that $\widetilde{S}$ is dense in
$\widetilde{Y}$
and
such that 
for any $\tilde{s} \in \widetilde{S}$,
each irreducible component
of $\pr_{A_1} \left( \pr_{B}^{-1} 
\left(
\tilde{s}_{\overline{K}}
\right) \cap X' \right) 
\subset
A_1$
has dense small points,
where $\tilde{s}_{\overline{K}} \in Y
\left(
\overline{K}
\right)$
is the point corresponding to $s$ via the natural map
$\widetilde{Y} (k) \hookrightarrow Y
\left(
\overline{K}
\right)$.

Note that
a general fiber of $\pr_{B}|_{X'} : X' \to Y$
has dimension $\dim (A_1)$.
Indeed, if this is not the case,
then
$\pr_{B}^{-1} 
\left(
\tilde{s}_{\overline{K}}
\right) \cap X'$ is an effective divisor of 
$\left\{
\tilde{s}_{\overline{K}}
\right\} \times A_1$
for any $s \in S$,
and hence
$\pr_{A_1} \left( \pr_{B}^{-1} 
\left(
\tilde{s}_{\overline{K}}
\right) \cap X' \right)$
is an effective divisor on $A_1$
each of which irreducible components has dense small points. 
However, 
since $A_1$ is nowhere degenerate and has trivial $\overline{K}/k$-trace,
this contradicts Step~1 above.
Thus a general fiber of $\pr_{B}|_{X'} : X' \to Y$
has dimension $\dim (A_1)$.

Since $X' \subsetneq B \times A_1$,
it follows that
$\widetilde{Y} \neq \widetilde{B}$.
Since 
$X'$ has codimension $1$,
this means $X' = \widetilde{Y} \times_{\Spec (k)} A_1$.
This shows that $X'$ is a special subvariety.

\smallskip

\emph{Step~3. For an arbitrary $A$.}
Let $\mathfrak{m}$ be the maximal nowhere degenerate abelian subvariety
of $A$.
Let $A \to A' := A/\mathfrak{m}$ be the quotient homomorphism.
By the 
Poincar\'e
complete reducibility theorem,
there exists an isogeny
$\varphi : A \to A' \times \mathfrak{m}$.
We set $\phi : A \to \mathfrak{m}$ to be the composite
$\pr_{\mathfrak{m}} \circ \varphi$,
where $\pr_{\mathfrak{m}} :  A' \times \mathfrak{m} \to \mathfrak{m}$
is the canonical projection.

We claim that $\phi (X)$ is a special subvariety of $\mathfrak{m}$.
If $\phi (X) = \mathfrak{m}$,
then it is special trivially.
Therefore, suppose that
$\phi (X) \subsetneq \mathfrak{m}$.
Then, since $X$ has codimension $1$ in $A$, so does 
$\phi (X)$ in $ \mathfrak{m}$.
Since $X$ has dense small points, so does $\phi (X)$
(cf. \cite[Lemma~2.1]{yamaki5}).
By Step~2 above,
it follows that $\phi (X)$ is special.

By assumption, furthermore, $X$ has dense small points.
By \cite[Theorem~7.21]{yamaki6}, it follows that
$X$ is a special subvariety of $A$.
This completes the proof of the theorem.
\QED

We give a couple of remarks.
The following proposition is shown in the Step~1 of the proof of Theorem~\ref{thm:divisorialcase}.

\begin{Proposition} \label{prop:nondegtrivialtrace}
Let $A$ be a nowhere degenerate abelian variety
over $\overline{K}$
with trivial $\overline{K}/k$-trace
and let $X$ be a closed subvariety of $A$ of codimension $1$.
Then $X$ does not have dense small points.
\end{Proposition}

\begin{Remark} \label{rem:main1-nondegtrivialtrace}
By Proposition~\ref{prop:dense-height0},
it follows from
Proposition~\ref{prop:nondegtrivialtrace}
that a non-trivial effective divisor
on a nowhere degenerate abelian variety over $\overline{K}$
with trivial $\overline{K}/k$-trace
has positive canonical height.
That is a generalization of Proposition~\ref{prop:main1}.
\end{Remark}

\begin{Remark} \label{rem:divisor-nontorsion}
Let $A$ be a nowhere degenerate abelian variety
over $\overline{K}$ with trivial $\overline{K}/k$-trace.
It follows from
Proposition~\ref{prop:nondegtrivialtrace}
that
no non-zero effective divisor
is a torsion subvariety of $A$,
but
we can show this directly, 
not via this proposition, in fact.
To prove that by contradiction,
suppose that there exists a torsion subvariety of codimension
$1$.
Then there exists an abelian subvariety $G$ of codimension $1$.
Take the quotient $A/G$.
Then by Lemma~\ref{lem:subandquotient},
it is nowhere degenerate
and has trivial $\overline{K}/k$-trace.
On the other hand,
since $A/G$ is a nowhere degenerate elliptic curve,
it follows from the well-known fact 
that the moduli space of
elliptic curves is affine that
$A/G$ is a constant abelian variety.
That is a contradiction.
\end{Remark}

Next, we consider the case where the subvariety has dimension $1$.
We show one more lemma:

\begin{Lemma} \label{lem:addedfinially}
Let $A$ be an abelian variety over $\overline{K}$
and let $X_1$ and $X_2$ be closed subvarieties of $A$ with dense small points.
Then 
$X_1 + X_2 := \{ x_1 + x_2 \mid x_1 \in X_1 , x_2 \in X_2 \}$
has dense small points.
\end{Lemma}

\Proof
Note that
$X_1 + X_2$ is the image 
of $X_1 \times X_2 \subset A \times A$ 
by the addition homomorphism
$\alpha : A \times A \to A$.
By \cite[Lemma~2.4]{yamaki5},
$X_1 \times X_2$ has dense small points, and by \cite[Lemma~2.1]{yamaki5},
it follows that $X_1 + X_2 = \alpha ( X_1 \times X_2 )$ 
has dense small points.
\QED

Now, we show the dimension $1$ case.

\begin{Theorem} \label{thm:curvecase}
Let $A$ be an abelian variety over $\overline{K}$
and
let $X$ be a closed subvariety of $A$ of dimension $1$.
If $X$ has dense small points, then it is a special subvariety.
\end{Theorem}

\Proof
\emph{Step~1. 
In this step,
we make an additional assumption that the $\overline{K}/k$-trace
of $A$ is trivial, under which
we prove 
that 
if $X$ has dense small points, then it is a torsion subvariety.}
Let $A \left( \overline{K} \right)_{tor}$ denote the set of torsion
points of $A \left( \overline{K} \right)$.
Consider
\[
\mathfrak{S} := 
\left\{
(\tau , A_{\tau})
\left|
\ 
\text{$\tau \in 
A \left(
\overline{K}
\right)_{tor}$,
$A_\tau$ is an abelian subvariety of $A$ with $X - \tau \subset A_{\tau}$}
\right.
\right\}
.
\]
We take $( \tau_{1} , A_{\tau_1}) \in \mathfrak{S}$
such that
\[
\dim \left( A_{\tau_1} \right)
=
\min
\left\{ 
\dim ( A_{\tau} )
\mid
\left( \tau , A_{\tau} \right) \in \mathfrak{S}
\right\}
.
\]
To ease notation,
we set $A_1 := A_{\tau_1}$ and $X_1 := X - \tau_1$.
Remark that $X_1 \subset A_1$, 
and that $A_1$ has trivial
$\overline{K}/k$-trace
by Lemma~\ref{lem:subandquotient}~(1).

We set $Y_0 := \{ 0 \} \subset A_{1}$ 
and $Y_1 := X_1 - X_1$.
For each integer $m >1 $, we define 
$Y_m$ to be the sum of $m$ copies of $Y_1$ in $A_1$:
\begin{align*}
Y_{m} := 
\underbrace{Y_1 + \cdots + Y_1}_{m} .
\end{align*}
Then,
since the sequence $(Y_m)_{m \in \ZZ_{\geq 0}}$ 
is an ascending sequence of closed 
(reduced irreducible) subvarieties
of $A_1$, this sequence is stable for large $m$,
that is, there exists a positive integer $N$
such that
$Y_{N-1} \subsetneq Y_{N} = Y_{N + 1} = Y_{N + 2} = \cdots$.
(Since $Y_0 \subsetneq Y_1$, we remark $N \geq 1$.)
Then we note that $Y_N$ is an abelian subvariety of $A_1$.
Indeed, $Y_N \neq \emptyset$,
and we see
$
Y_N - Y_N \subset Y_{2N} = Y_N
$.

Now, suppose that $X$ has dense small points.
We would like to show that $X_1 = A_1$.
Note that $X_1 = X - \tau_1$
also has dense small points.

\begin{Claim} \label{cl:Y_N=A_1}
We have $Y_N = A_1$.
\end{Claim}

\Proof
Let
$\psi : A_1 \to A_{1} / Y_N$
be the quotient homomorphism.
Since $A_1$ has trivial $\overline{K}/k$-trace,
so does
$A_{1} / Y_N$ by Lemma~\ref{lem:subandquotient}~(1).
Fix a point $x_1 \in X_1 \left( \overline{K} \right)$.
Then $X_1 - x_1 \subset Y_1 \subset Y_N$.
It follows that
$\psi ( X_1 - x_1 ) = \{ 0 \}$,
and hence
$\psi (X_1) = \{ \psi (x_1) \}$.
Since $X_1$ has dense small points, so does $\{ \psi (x_1) \}$
by \cite[Lemma~2.1]{yamaki5},
which means that $\psi (x_1)$ is a torsion point.
Since a surjective morphism of abelian varieties over an algebraically closed
field induces a surjective morphism between the torsion points,
there exists a $\tau' \in A_1 \left( \overline{K}
\right)_{tor}$ such that $\psi ( \tau' ) = \psi (x_1)$.
Then $\psi ( X_1 - \tau' ) = \{ 0 \}$, which means
$Y_N \supset X_1 - \tau' = X - (\tau_1 + \tau')$.
We note
that $\tau_1 + \tau'$ is a torsion point.
By the minimality of $\dim (A_1)$ among
$\left\{ 
\dim ( A_{\tau} )
\mid
\left( \tau , A_{\tau} \right) \in \mathfrak{S}
\right\}$,
we have $\dim (Y_N) \geq \dim (A_1)$.
Since $Y_N \subset A_1$, this proves $Y_N = A_1$.
\QED

\begin{Claim} \label{cl:N=1}
We have $N = 1$.
\end{Claim}

\Proof
We argue by contradiction.
Suppose that $N \geq 2$.
Then $Y_1 \subset Y_{N-1} \subsetneq A_1$.
Since  $X_1$ has dense small points,
so does
$Y_{N-1}$ by Lemma~\ref{lem:addedfinially}.

Note that $Y_{N-1}$ is not an abelian subvariety.
Indeed, if this is not the case,
then $Y_{N-1} + Y_{N-1} \subset Y_{N-1}$,
so that
$Y_{N} = Y_{N-1} + Y_1 \subset Y_{N-1} + Y_{N-1} \subset Y_{N-1}$,
which contradicts the choice of $N$.

Since $A_1 = Y_N = Y_{N-1} + Y_1$ 
(cf. Claim~\ref{cl:Y_N=A_1}) and $1 \leq \dim Y_1 \leq 2$,
we have
$1 \leq \codim \left( Y_{N-1} , A_1 \right) \leq 2$.
In fact, 
we see that $\codim \left( Y_{N-1} , A_1 \right) = 2$.
To show that by contradiction,
suppose that $\codim \left( Y_{N-1} , A_1 \right) = 1$.
Since $Y_{N-1}$ has dense small points,
$Y_{N - 1}$ is then a special subvariety by Theorem~\ref{thm:divisorialcase}.
Since the $\overline{K}/k$-trace of $A_1$
is trivial by the assumption in this step
and since $0 \in Y_{N-1}$,
it follows that $Y_{N-1}$ is an abelian subvariety.
However, this contradicts what we have noted above.
Thus $\codim \left( Y_{N-1} , A_1 \right) = 2$.

Set $Z := Y_{N-1} + X_1$.
Then
by Lemma~\ref{lem:addedfinially},
it is a closed subvariety of $A_1$ with dense small points.
Furthermore,
since $\codim \left( Y_{N-1} , A_1 \right) = 2$
and $Z -X_1 = Y_N = A_1$,
we have $\codim (Z ,A_1) = 1$.
Again by Theorem~\ref{thm:divisorialcase}
together with the fact that the $\overline{K}/k$-trace
of $A_1$ is trivial,
it turns out that $Z$ is a torsion subvariety.
Thus there exist an abelian subvariety $G$ of $A_1$
and a torsion point $\tau' \in A_1 \left(
\overline{K} \right)$
such that $Z = G + \tau'$.
Since $0 \in Y_{N-1}$, we have $X_1 \subset Y_{N-1} + X_1 = Z$.
It follows that $X - (\tau' + \tau_1) =
X_{1} - \tau' \subset G$.
By the minimality of $\dim (A_1)$, this implies that $\dim (G) \geq \dim (A_1)$.
Since $G \subset A_1$, we obtain $G = A_1$,
and thus $Z = A_1$.
However, this contradicts $\codim (Z ,A_1) = 1$.
Thus we obtain $N = 1$.
\QED

By Claim~\ref{cl:Y_N=A_1}
and Claim~\ref{cl:N=1}, we have $X_1 - X_1 = A_1$.
Since $\dim (X_1) = 1$, we have $1 \leq \dim (A_1) \leq 2$.
\begin{Claim} \label{cl:dimA_1=1}
We have $\dim (A_1) = 1$.
\end{Claim}

\Proof
We argue
by contradiction.
Suppose that $\dim (A_1) \neq 1$, 
or equivalently,
$\dim (A_1) = 2$.
Then $X_1$ is a divisor on $A_1$.
Since $X_1$ has dense small points,
it follows from
Theorem~\ref{thm:divisorialcase}
that $X_1$ is a torsion subvariety,
where we note that the $\overline{K}/k$-trace of $A_1$ 
is trivial.
Since $\dim (X_1) = 1$,
this means that the translate of $X_1$ by some torsion point
is an abelian subvariety of dimension $1$,
which contradicts the definition of $A_1$.
\QED

Since $X_1 \subset A_1$,
we have $X_1 = A_1$ by Claim~\ref{cl:dimA_1=1}.
Thus $X = X_{1} + \tau_1 = A_{1} + \tau_1$,
which is a torsion subvariety of $A$.
This completes Step~1.

\medskip
\emph{Step~2. We consider the general case.}
Let $\mathfrak{t}$ be the image of 
the
$\overline{K}/k$-trace homomorphism (cf. \S~\ref{sect:NC}).
Let $\phi : A \to A / \mathfrak{t}$ be the quotient homomorphism
and set $X' := \phi (X)$.
Since $X$ has dense small points , so does $X'$
by \cite[Lemma~2.1]{yamaki5}.
Since $A / \mathfrak{t}$ has trivial $\overline{K}/k$-trace
(cf. \cite[Remark~5.4]{yamaki7}),
it follows from
Remark~\ref{rem:trivialcase}
and
Step~1 that there exist an abelian subvariety $G \subset A/\mathfrak{t}$ 
and a torsion point $\tau' \in A/\mathfrak{t}$
such that $X' = G + \tau'$.
Note $\dim (G) \leq 1$.
Since $\phi$ induces a surjective homomorphism between the torsion points,
we take a torsion point $\tau \in A \left( 
\overline{K} \right)$ with $\phi ( \tau ) = \tau'$.

Set $Y := X - \tau$
and $B := \phi^{-1} (G)$.
Then $Y \subset B$,
and
since $\phi : A \to A / \mathfrak{t}$ is smooth and has connected fibers,
$B$ is an abelian subvariety of $A$
with $\mathfrak{t} \subset B$.
Let $\mathfrak{n}$ be the maximal nowhere degenerate abelian 
subvariety of $B$
and let $\mathfrak{s}$ be the image of the 
$\overline{K}/k$-trace 
homomorphism
of $B$.
By the universality of the $\overline{K}/k$-trace
of $B$,
we have $\mathfrak{t} \subset \mathfrak{s}$.
(The other inclusion holds by the universality of
the
$\overline{K}/k$-trace
of $A$ also,
and thus
 they are equal in fact.)
Since then $\dim ( \mathfrak{n} / \mathfrak{s} )
\leq
\dim ( \mathfrak{n} / \mathfrak{t} ) 
\leq 
\dim (G)
\leq 1$,
the geometric Bogomolov conjecture holds for $\mathfrak{n} / \mathfrak{s}$
(cf. Remark~\ref{rem:trivialcase}).
By \cite[Theorem~5.5]{yamaki7},
it follows that the conjecture holds for $B$.
Since $Y$ has dense small points,
$Y$ is therefore a special subvariety.
Thus $X = Y + \tau$ is also a special subvariety,
which completes the proof of the theorem.
\QED

\subsection{Bogomolov conjecture for curves} \label{subsect:GBCforC}

In this subsection,
let $C$
be a smooth projective curve of genus $g \geq 2$
over $\overline{K}$.
Let $J_C$ be the Jacobian variety of $C$.
Fix a divisor $D$ of degree $1$ on $C$
and let
$
\jmath_{D} : C \hookrightarrow 
J_C
$
be the embedding
defined
by $\jmath_D (x) := x - D$.
For any $P \in J_{C} \left( \overline{K} \right)$
and
any $\epsilon \geq 0$,
set
\[
B_{C} ( P , \epsilon )
:=
\left\{ 
\left. x \in C \left(
\overline{K}
\right) 
\ 
\right|
|| \jmath_D(x) - P ||_{NT} \leq \epsilon
\right\}
,
\]
where $|| \cdot ||_{NT}$ is the semi-norm arising from the N\'eron--Tate
height on $J_{C}$
associated to a symmetric theta divisor.

Recall that $C$ is said to  be \emph{isotrivial}
if there exists a projective curve $\widetilde{C}$ over $k$
with $C \cong \widetilde{C} \otimes_k \overline{K}$.
As is remarked in \cite[\S~8]{yamaki6},
we obtain the following theorem
as a consequence of 
Theorem~\ref{thm:curvecase}.

\begin{Theorem} [Theorem~\ref{thm:main1intro},
Bogomolov conjecture for curves
over any function field]
\label{thm:BCFC}
Assume that $C$ is non-isotrivial.
Then, for any $P \in J_{C} \left(
\overline{K} \right)$, there exists an $r > 0$ such that
$B_{C} ( P , r )$ is
finite.
\end{Theorem}

Thus we conclude that the Bogomolov conjecture
for curves over any function field holds.

\begin{Remark}
In the above theorem, we consider the N\'eron--Tate height associated to
a symmetric theta divisor,
but this theorem also holds with respect to
the canonical height associated to any even ample line bundle
because of
\cite[Lemma~2.1]{yamaki5}.
\end{Remark}

\section{Consequences of the main results} \label{sect:RGBC}

\subsection{Geometric Bogomolov conjecture
for abelian varieties} \label{subsect:GBC}

In this subsection,
we make a contribution to the geometric Bogomolov conjecture
(Conjecture~\ref{conj:GBCforAVint}),
which claims that
a closed subvariety $X$ of $A$ has dense small points
if only if $X$ is a special subvariety.
Since we know the ``if'' part holds (cf. \cite{yamaki5}),
the ``only if'' part is the problem.

We begin with some background of the geometric Bogomolov conjecture.
The version 
of the conjecture 
over number fields was first established by Zhang
in \cite{zhang2},
which claims that for a closed subvariety
$X$
of an abelian variety over a number field,
$X$ has dense small points if and only if $X$ is a torsion subvariety.
Zhang's theorem generalizes  Ullmo's theorem,
which is the case of a curve in its Jacobian,
and its proof is inspired by the proof of Ullmo.
Further,
Moriwaki generalized the result to the case
where $K$ is a finitely generated field
over $\QQ$,
with respect to certain arithmetic heights.
Note that Moriwaki's arithmetic heights are not the classical
heights over function fields.

Over function fields
with respect to the classical heights, no analogous results were known
for a while
even after Zhang's theorem,
but in 2007,
Gubler proved the following result.
Here $A$ is \emph{totally degenerate} at $v \in M_{\overline{K}}$
if $A$ has torus reduction at $v$.

\begin{citeTheorem} [Theorem~1.1 of \cite{gubler2}]
Let $K$ be a function field.
Let $A$ be an abelian variety over $\overline{K}$.
Assume that $A$ is totally degenerate at some $v \in M_{\overline{K}}$.
Let $X$ be a closed subvariety of $A$.
Then $X$ has dense small points if and only if $X$ is a
torsion subvariety.
\end{citeTheorem}

The geometric Bogomolov conjecture is inspired by
Gubler's theorem.
Remark that his theorem is equivalent to the geometric
Bogomolov conjecture for an abelian variety which is
totally degenerate at some place
because then a special subvariety is a torsion subvariety.

Although the geometric Bogomolov conjecture for abelian varieties
is still an open problem,
there are partial solutions.
Recently, in \cite{yamaki6} and \cite{yamaki7},
we obtain results which
generalize Gubler's theorem.
Here we would like to recall them.
Let $A$ be an abelian variety over $\overline{K}$
and let
$\mathfrak{m}$ be 
the maximal nowhere degenerate abelian subvariety of $A$.
We remark that if $A$ is totally degenerate at some place,
then $\mathfrak{m} = 0$
(cf. \cite[Introduction]{yamaki6}), 
but the converse does not holds in general.
Let $\mathfrak{t}$ be the image of the trace homomorphism
$\Tr_A : \widetilde{A}^{\overline{K}/k} \otimes_k \overline{K}
\to A$. 
Since $\widetilde{A}^{\overline{K}/k} \otimes_k \overline{K}$
is nowhere degenerate,
$\mathfrak{t} \subset \mathfrak{m}$
(cf. \cite[Lemma~7.8~(2)]{yamaki6}),
and $\mathfrak{m} /\mathfrak{t} $
is a nowhere degenerate abelian variety with trivial $\overline{K}/k$-trace
(cf. \cite[Lemma~7.8~(2)]{yamaki6} and \cite[Remark~5.4]{yamaki7}).
Following the work \cite{yamaki6}, we obtain 
in \cite{yamaki7} the following result.

\begin{citeTheorem} [cf. Theorem~1.5 of \cite{yamaki7}]
\label{thm:yamaki72int}
Let $A$ be an abelian variety over $\overline{K}$.
Then the following are equivalent to each other:
\begin{enumerate}
\renewcommand{\labelenumi}{(\alph{enumi})}
\item
The geometric Bogomolov conjecture holds for $A$;
\item
The geometric Bogomolov conjecture holds for $\mathfrak{m}/\mathfrak{t}$.
\end{enumerate}
\end{citeTheorem}

Theorem~\ref{thm:yamaki72int}
shows that
the geometric Bogomolov
conjecture for abelian varieties 
is equivalent to
that for nowhere degenerate abelian varieties
with trivial $\overline{K}/k$-trace. 
In particular, Theorem~\ref{thm:yamaki72int}
proves
that the conjecture holds for $A$ with 
$\dim (\mathfrak{m} /  \mathfrak{t}) = 0$
(cf. \cite[Theorem~1.4]{yamaki6}).
Since
$\mathfrak{m} =  \mathfrak{t} = 0$
for a totally degenerate abelian variety $A$,
this 
generalizes Gubler's theorem.

Our main results 
of this paper give a further generalization of the above results.
First,
as a consequence of Theorems~\ref{thm:divisorialcase}
and \ref{thm:curvecase},
we obtain the following result
(see also Remark~\ref{rem:trivialcase}).

\begin{Corollary} \label{cor:Aleq3}
Let $A$ be an abelian variety over $\overline{K}$
with $\dim (A) \leq 3$.
Then
the geometric Bogomolov conjecture holds for $A$.
\end{Corollary}

Then,
by Theorem~\ref{thm:yamaki72int},
we generalize Corollary~\ref{cor:Aleq3}
as follows.

\begin{Corollary} \label{cor:mleqt+3}
Let $A$ be an abelian variety over $\overline{K}$.
Let $\mathfrak{m}$ be the maximal nowhere degenerate abelian variety
and let $\mathfrak{t}$ be the image of the 
trace homomorphism.
Assume that $\dim (\mathfrak{m}) 
\leq
\dim (\mathfrak{t}) 
+ 3$.
Then the geometric Bogomolov conjecture holds for $A$.
\end{Corollary}

\begin{Remark} \label{rem:trace-finite}
Since the trace homomorphism is finite,
we have $\dim (\mathfrak{t}) = \dim \left( \widetilde{A}^{\overline{K}/k} \right)$.
\end{Remark}

The geometric Bogomolov conjecture
remains open even now for a higher dimensional abelian variety,
but we can see that the conjecture holds for some kind of
them.
For example,
let $A$ be an abelian variety 
of dimension $4$
over $\overline{K}$.
Then we have proved that
the geometric Bogomolov conjecture holds for it
unless $A$ is a nowhere degenerate abelian variety with 
trivial $\overline{K}/k$-trace.

\subsection{Manin--Mumford conjecture and the Bogomolov conjecture}
\label{subsect:scanlontheorem}

In this subsection,
we remark that our main results 
give an alternative proof of the
Manin--Mumford conjecture
in positive characteristic
in a special setting,
by mentioning
the relationship between 
this conjecture
and the geometric Bogomolov conjecture.

Let $F$ be an algebraically closed field.
Originally, the Manin--Mumford conjecture claims that
under the assumption of $\ch (F) = 0$,
a smooth projective curve of genus at least $2$
over $F$
embedded in its Jacobian
should have only finitely many torsion points of the Jacobian.

This was first proved by Raynaud in \cite{raynaud1}
in the following generalized form.
Here, we say
that a closed subvariety $X$ of an abelian variety $A$
has \emph{dense torsion points}
if 
$ X \cap A ( F )_{tor}$ is dense in $X$.

\begin{citeTheorem} \label{thm:raynaud1}
Let $X$ be a closed subvariety.
Assume that $\dim (X) = 1$.
If $X$ has dense torsion points,
then $X$ is a torsion subvariety.
\end{citeTheorem}

Further,
Raynaud proved
in \cite{raynaud2} 
that
the above statements holds for any dimensional $X$.
We should remark that Moriwaki proved in \cite{moriwaki5}
that 
Raynaud's theorem is a consequence of his theorem
on the 
arithmetic Bogomolov conjecture over a 
field
finitely
generated over $\QQ$.

Now, we consider what happens when $\ch (F) = p >0$.
In this case,
there exists an obvious counterexample:
If $F$ is an algebraic closure of a finite field
and $X$ is any closed subvariety of $A$,
then any $F$-point of $X$ is a torsion point,
and thus $X$ always has dense torsion points.
However,
up to such influence of subvarieties which can be
defined over finite fields,
one may expect 
that a similar statement should also hold in positive characteristics.
We have the following precise result.

\begin{citeTheorem} \label{thm:MM}
Let $F$ be an algebraically closed field
with $\ch (F) >0$
and let $k$ be the algebraic closure in $F$ 
of the prime field of $F$.
Let $A$ be an abelian variety over $F$
and 
let $X$ be a closed subvariety of $A$.
If $X$ has dense torsion points,
then there exist an abelian subvariety $G$ of $A$,
a closed subvariety $\widetilde{Y}$ of $\widetilde{A}^{F/k}$,
and a torsion point $\tau \in A \left( F \right)$ 
such that $X = G + \Tr^{F/k}_{A} \left( \widetilde{Y}
\otimes_k F \right) + \tau$,
where $\left( 
\widetilde{A}^{F/k} , \Tr^{F/k}_{A}
\right)$
is the $F/k$-trace of $A$.
\end{citeTheorem}

This theorem is due to the following authors:
In 2001,
Scanlon gave a sketch of the model-theoretic proof of this theorem
(\cite{scanlon0}).
In 2004,
Pink and Roessler gave an algebro-geometric proof
(\cite{PR}).
In 2005, Scanlon gave a detailed model-theoretic proof
in \cite{scanlon1}
based on the argument in \cite{scanlon0}.
Note that
in those papers, they prove a generalized version 
for semiabelian varieties $A$, in fact.

Here, we mention that Theorem~\ref{thm:MM} can be deduced from the geometric
Bogomolov conjecture for $A$, as Moriwaki did in the case of characteristic $0$.
Let $F$, $k$, $A$,
and $X$ be as in Theorem~\ref{thm:MM}.
Then there exist $t_1 , \ldots , t_n \in F$
such that $A$ and $X$ can be defined over $K := k ( t_1 , \ldots , t_n)$,
that is,
there exist an abelian variety $A_0$ over $K$
and a closed subvariety $X_0$ of $A_0$ such that
$A = A_0 \otimes_K F$ and $X = X_0 \otimes_K F$.
Further,
there exists a normal projective variety $\mathfrak{B}$ over $k$
with function field $K$.
Let $\overline{K}$ be the algebraic closure of $K$ in $F$
and set $A_{\overline{K}} := A_0 \otimes_K \overline{K}$
and $X_{\overline{K}} := X_0 \otimes_K \overline{K}$.
Let $\mathcal{H}$ be an ample line bundle on $\mathfrak{B}$.
Then we have a notion of height
over $K$, and we can consider the canonical height
on 
$A_{\overline{K}}$
associated to an even ample line bundle.
Suppose that $X$ has dense torsion points.
Then $X_{\overline{K}}$ has dense torsion points
and hence has dense small points.
Then if we assume the geometric Bogomolov conjecture,
it follows that $X_{\overline{K}}$ should be a special subvariety,
and this should show the conclusion of Theorem~\ref{thm:MM}.

Since 
the geometric Bogomolov conjecture 
is still open, the above argument 
does not give a new proof of Theorem~\ref{thm:MM}.
However,
we can actually deduce
the theorem in the following special cases:
\begin{enumerate}
\item
$\dim (X) = 1$ or $\codim (X , A) = 1$
(from Theorems~\ref{thm:curvecase}, \ref{thm:divisorialcase});
\item
$\dim (A) \leq 3$
(from Corollary~\ref{cor:Aleq3}).
\end{enumerate}
In particular, we have recovered the positive characteristic version
of Theorem~\ref{thm:raynaud1}.


\renewcommand{\thesection}{Appendix} 
\renewcommand{\theTheorem}{A.\arabic{Theorem}}
\renewcommand{\theClaim}{A.\arabic{Theorem}.\arabic{Claim}}
\renewcommand{\theequation}{A.\arabic{Theorem}.\arabic{Claim}}
\renewcommand{\theProposition}
{A.\arabic{Theorem}.\arabic{Proposition}}
\renewcommand{\theLemma}{A.\arabic{Theorem}.\arabic{Lemma}}
\setcounter{section}{0}
\renewcommand{\thesubsection}{A.\arabic{subsection}}




\small{

}

\end{document}